\documentclass[a4paper,10pt,leqno]{amsart}
\usepackage{enumerate,amssymb}
\usepackage[arrow,curve,matrix,tips,2cell]{xy}
  \SelectTips{eu}{10} \UseTips
  \UseAllTwocells

\DeclareMathAlphabet{\matheurm}{U}{eur}{m}{n}

\DeclareMathOperator{\CAT}{CAT}
\DeclareMathOperator{\Cay}{Cay}
\DeclareMathOperator{\cd}{cd}
\DeclareMathOperator{\cent}{cent}

\DeclareMathOperator{\HNN}{HNN}
\DeclareMathOperator{\id}{id}

\DeclareMathOperator{\Nil}{Nil}
\DeclareMathOperator{\Out}{Out}

\DeclareMathOperator{\res}{res}
\DeclareMathOperator{\rk}{rk}

\DeclareMathOperator{\Sol}{Sol}

\DeclareMathOperator{\vcd}{vcd}

  \newcommand{\IC}{\mathbb{C}}

  \newcommand{\IH}{\mathbb{H}}

  \newcommand{\IQ}{\mathbb{Q}}
  \newcommand{\IR}{\mathbb{R}}

  \newcommand{\IZ}{\mathbb{Z}}

  \newcommand{\calt}{\mathcal{T}}

\newcommand{\eub}[1]{\underline{E}#1}

\newcommand{\xycomsquare}[8]                      
{\xymatrix{#1 \ar[r]^{#2} \ar[d]^{#4} &
#3 \ar[d]^{#5}  \\
#6\ar[r]^{#7} &
#8
}
}

\newcommand{\xycomsquareminus}[8]                      
{\xymatrix{#1 \ar[r]^-{#2} \ar[d]^-{#4} &
#3 \ar[d]^-{#5}  \\
#6\ar[r]^-{#7} &
#8
}
}

\newcounter{commentcounter}

\theoremstyle{plain}
\newtheorem{theorem}{Theorem}[section]

\newtheorem{lemma}[theorem]{Lemma}

\newtheorem{conjecture}[theorem]{Conjecture}

\theoremstyle{definition}
\newtheorem{definition}[theorem]{Definition}

\newtheorem{example}[theorem]{Example}

\newtheorem{remark}[theorem]{Remark}

\theoremstyle{remark}

\makeatletter\let\c@equation=\c@theorem\makeatother

\newcommand{\action}{\curvearrowright}

\title{Survey on geometric group theory}
              \author{Wolfgang L\"uck}
      \address{Westf\"alische Wilhelms-Universit\"at M\"unster\\
               Mathematisches Institut\\
               Einsteinstr.~62,
               D-48149 M\"unster, Germany}
                \email{lueck@math.uni-muenster.de}
      \urladdr{http://www.math.uni-muenster.de/u/lueck}
         \date{September 2008}
     \keywords{geometric group theory, measurable group theory, quasiisometry,
hyperbolic groups, classifying spaces}
    \subjclass[2000]{20-02, 20F65, 20F67, 20F69, 57-02, 57M60, 57P10}

\begin{document}


\begin{abstract}
This article is a survey article on geometric group theory from the point
of view of a non-expert who likes geometric group theory and uses it in his own research.
\end{abstract}

\maketitle


\typeout{--------------------   Section 0: Introduction --------------------------}

\setcounter{section}{-1}
\section{Introduction}
\label{sec:Introduction}

This survey article on geometric group theory is written by a non-expert 
who likes geometric group theory 
and uses it in his own research. It is meant as a service
for people who want to receive an impression and read an introduction about the topic
and possibly  will later pass to more elaborate and specialized survey articles or to
actual research articles. There will be no proofs. Except
for Theorem~\ref{the:hyperbolic_groups_and_manifolds} all results have already appeared
in the literature.

There is to the author's knowledge no obvious definition what 
geometric group theory really is. At any rate the basic idea 
is to pass from a finitely generated group to the geometry
underlying its Cayley graph with the word metric. It turns out that
only the large scale geometry is really an invariant of the group itself
but that this large scale or coarse geometry carries a lot of information.
This leads also to a surprising and intriguing variety of new results and structural 
insights about groups. 

A possible explanation for this may be
that humans have a better intuition when they think in geometric terms.
Moreover, it is helpful to understand  groups in the way 
as they have appeared naturally in mathematics,
namely, as  groups of symmetries. In other words, basic information about
a group can be obtained by studying its actions on nice spaces.

The personal interest of the author comes from questions of the type
whether a group satisfies the conjectures due to Baum-Connes, Borel, 
Farrell-Jones, Kaplansky, Novikov, Hopf, Singer or yields a positive 
answer to Atiyah's question
on $L^2$-Betti numbers. They are all of the kind that one wants to know
whether for a given group $G$ its group ring $RG$, its reduced group $C^*$-algebra
$C^*_r(G)$, or an aspherical closed manifold with $G$ as fundamental group satisfy
certain algebraic or geometric properties concerning their structure
as rings or $C^*$-algebras, their $K$- or $L$-theory, rigidity properties
or the spectrum of the Laplace operator of the universal covering. A priori
these problems do not seem to be related to questions about the geometry of the group.
However, most of the proofs for certain  classes of groups contain an important part, where one uses
certain geometric properties of the groups, very often properties such
as being negatively or  non-positively curved in some metric sense. For instance, 
there is the, on first sight purely ring theoretic, 
conjecture that for a torsionfree group $G$ and an integral domain $R$ 
the group ring $RG$ contains no idempotents except $0$ and $1$. It is 
surprising that a proof of it can be given for certain rather large classes of groups by 
exploiting their geometry, and no algebraic proof is known in these cases.

The author has done his best to sort out interesting problems and results and to include
the relevant references, and apologizes if an important aspect or reference
is missing, it was left out because of ignorance, not on purpose.

The work was financially supported by the
Sonderforschungsbereich 478 \--- Geometrische Strukturen in der
Mathematik \--- and the Max-Planck-Forschungspreis and the Leibniz-Preis of the
author. The author wishes to thank Tom Church, Jan Essert, Ralf Gramlich, Clara L\"oh,
Sayed Roushon, Roman Sauer and Yehuda Shalom for their useful comments and 
in particular the two referees for  their very valuable detailed reports.

The paper is organized as follows:
\\[2mm]
\begin{tabular}{ll}%
\ref{sec:classical_examples}. 
& Classical examples 
\\%
\ref{sec:Basics_about_quasiisometry}. 
& Basics about quasiisometry
\\%
\ref{sec:Properties_and_invariants_of_groups_invariant_under_quasiisometry}. 
& Properties and invariants of groups invariant under quasiisometry
\\%
\ref{sec:Rigidity}. 
& Rigidity
\\%
\ref{sec:Hyperbolic_spaces_and:cat(kappa)-spaces}. 
& Hyperbolic  spaces and CAT($\kappa$)-spaces
\\%
\ref{sec:The_boundary_of_a-hyperbolic_spaces}. 
& The boundary of a hyperbolic  space
\\%
\ref{sec:Hyperbolic_groups}. 
& Hyperbolic groups
\\%
\ref{sec:CAT(0)-groups}. 
& CAT(0)-groups
\\%
\ref{sec:Classifying_spaces_for_proper_actions}. 
& Classifying spaces for proper actions
\\%
\ref{sec:Measurable_group_theory}. 
& Measurable group theory
\\%
\ref{sec:Some_open_problems}. 
& Some open problems
\\
& References
\end{tabular}


\typeout{--------------------   Section 1 ---------------------------------------}

\section{Classical examples}
\label{sec:classical_examples}

A classical example of geometric methods used in group theory is the topological
proof of Schreier's theorem.

\begin{theorem}[Schreier's Theorem]
\label{the:Schreier}
Let $G$ be a free group and $H \subseteq G$ be a subgroup. 
Then $H$ is free. If the rank $\rk(G)$ and the index $[G:H]$ are finite, then the rank of
$H$ is finite and satisfies
$$\rk(H) = [G:H] \cdot \bigl(\rk(G) -1\bigr) + 1.$$ 
\end{theorem}
\begin{proof}
Let $G$ be a free group on the set $S$. Take the wedge 
$X = \bigvee_{S} S^1$ of  circles, one copy for each element in $S$.
This is a $1$-dimensional $CW$-complex with $\pi_1(X) \cong G$ by
the Seifert-van Kampen Theorem. Let $p \colon \overline{X}
\to X$ be the covering associated to $H \subseteq G = \pi_1(X)$. We have
$\pi_1(\overline{X}) \cong H$. Since $X$ is a $1$-dimensional $CW$-complex,
$\overline{X}$ is a $1$-dimensional $CW$-complex. If $T \subseteq \overline{X}$
is a maximal tree, then $\overline{X}$ is homotopy equivalent to $
\overline{X}/T = \bigvee_{\overline{S}} S^1$
for some set $\overline{S}$. By the Seifert-van Kampen Theorem $H \cong \pi_1(\overline{X})$
is the free group generated by the set $\overline{S}$.

Suppose that $\rk(G)$ and $[G : H]$ are finite. Since $|S| = \rk(G)$,
the $CW$-complex~$X$ is compact. Since  $[G : H]$ is finite, the
$CW$-complex~$\overline{X}$ 
and hence $\overline{X}/T$ are compact. Hence 
$\rk(H) = |\overline{S}|$ is finite. We obtain for the Euler characteristics
$$1 - |\overline{S}| = \chi(\overline{X}) = [G:H] \cdot \chi(X) = 
[G:H] \cdot \left(1 - |S|\right).$$
Since $|S| = \rk(G)$ and $|\overline{S}| = \rk(H)$, the claim follows.
\end{proof}

Another example of this type is the topological proof of Kurosh's Theorem,
which can be found for instance 
in~\cite[Theorem~14 in~I.5 on page~56]{Serre(1980)}. The interpretation of
amalgamated products and $\HNN$-extensions in terms of topological spaces by
the Seifert-van Kampen Theorem or actions of groups on trees are in the same spirit
(see for instance~\cite{Baumslag(1993)},\cite{Cohen(1989)},
\cite{Dicks-Dunwoody(1989)}, \cite{Lyndon-Schupp(1977)}, 
\cite{Serre(1980)}).


\typeout{--------------------   Section 2 ---------------------------------------}

\section{Basics about quasiisometry}
\label{sec:Basics_about_quasiisometry}

A very important notion is the one of quasiisometry since it yields a bridge
between group theory and geometry by assigning to a finitely generated group 
a metric space (unique up to quasiisometry), namely, its Cayley graph with the word metric.
There are many good reasons for this passage, see 
for instance the discussion in~\cite[Item~0.3 on page 7 ff.]{Gromov(1993)}). 
At any rate this concept
has led to an interesting and overwhelming variety of new amazing results and applications
and to intriguing and stimulating activities.

\begin{definition} \label{def:quasiisometry}
Let $X_1 = (X_1,d_1)$ and $X_2 = (X_2,d_2)$ be two metric spaces. A map $f \colon X_1 \to X_2$ 
is called a \emph{quasiisometry} if there exist real numbers $\lambda, C > 0$ satisfying:

\begin{enumerate}

\item The inequality
$$\lambda^{-1} \cdot d_1(x,y) - C \le d_2\bigl(f(x),f(y)\bigr) \le \lambda \cdot  d_1(x,y) + C$$
holds for all $x,y \in X_1$;

\item For every $x_2$ in $X_2$ there exists $x_1 \in X_1$ with $d_2(f(x_1),x_2) < C$.

\end{enumerate}

We call $X_1$ and $X_2$ \emph{quasiisometric} if there is a quasiisometry
$X_1 \to X_2$.
\end{definition}

\begin{remark}[Quasiisometry is an equivalence relation]
If $f \colon X_1 \to X_2$ is a quasiisometry, then there exists a quasiisometry
$g \colon X_2 \to X_1$ such that both composites $g \circ f$ 
and $f \circ g$ have bounded distance from the identity map. The composite
of two quasiisometries is again a quasiisometry. Hence the notion of quasiisometry
is an equivalence relation on the class of metric spaces.
\end{remark}

\begin{definition}[Word-metric]
Let $G$ be a finitely generated group. Let $S$ be a finite 
set of generators.  The \emph{word metric} 
$$d_S \colon G \times G \to \IR$$
assigns to $(g,h)$ the minimum over all integers $n \ge 0$ such that
$g^{-1}h$ can be written as a product $s_1^{\epsilon_1}s_2^{\epsilon_2} \ldots s_n^{\epsilon_n}$
for elements $s_i \in S$ and $\epsilon_i \in \{\pm 1\}$.
\end{definition}

The metric $d_S$ depends on $S$. The main motivation for the notion of quasiisometry
is that the quasiisometry class of $(G,d_S)$ is independent  of the choice of 
$S$ by the following elementary lemma.

\begin{lemma} \label{lem:quasiisometry_invarince_word_metric}
Let $G$ be a finitely generated group.
Let $S_1$ and $S_2$ be two finite sets of generators. 
Then the identity $\id \colon (G,d_{S_1}) \to (G,d_{S_2})$ is a quasiisometry.
\end{lemma}
\begin{proof}
Choose $\lambda$ such that for all $s_1 \in S_1$ we have $d_{S_2}(s_1,e), d_{S_2}(s_1^{-1},e)
\le \lambda$
and for $s_2 \in S_2$ we have $d_{S_1}(s_2,e), d_{S_1}(s_2^{-1},e) \le \lambda$.
Take $C = 0$.
\end{proof}

\begin{definition}[Cayley graph] 
\label{def:Cayley_graph}
Let $G$ be a finitely generated group. Consider a finite set $S$ of generators.
The \emph{Cayley graph} $\Cay_S(G)$ is the graph whose set of vertices is $G$
and there is an  edge joining $g_1$ and $g_2$ if and only if
$g_1 = g_2s$ for some $s \in S$. 
\end{definition}

A \emph{geodesic} in a metric space $(X,d)$ is an isometric embedding $I \to X$,
where $I \subset \IR$ is an interval equipped with the metric induced 
from the standard metric on $\IR$. 

\begin{definition}[Geodesic space]
\label{def:geodesic_space}
A metric space $(X,d)$ is called a \emph{geodesic space}
if for two points $x,y \in X$ there is a geodesic $c \colon [0,d(x,y)] \to X$
with $c(0) = x$ and $c(d(x,y)) = y$.  
\end{definition}

Notice that we do not require the unique existence of a geodesic joining two given points.

\begin{remark}[Metric on the Cayley graph]
\label{rem:metric_on_the_Cayle_graph}
There is an obvious procedure to define a metric on $\Cay_S(G)$ such that
each edge is isometric to $[0,1]$ and such that the distance of two points
in $\Cay_S(G)$ is the infimum over the length over all piecewise linear paths joining
these two points. This metric restricted to $G$ is just the word metric $d_S$.
Obviously the inclusion $(G,d_S) \to \Cay_S(G)$ is a quasiisometry. In particular,
the quasiisometry class of the metric space $\Cay_S(G)$ is independent of $S$.

The Cayley graph allows to translate properties of a finitely
generated group to properties
of a geodesic metric space.
\end{remark}

\begin{lemma}[\v{S}varc-Milnor Lemma] \label{lem:Svarc-Milnor_Lemma}
Let $X$ be a geodesic space. 
Suppose that $G$ acts properly, cocompactly and isometrically 
on $X$. Choose a base point $x \in X$. Then the map
$$f \colon G \to X, \quad g \mapsto gx$$
is a quasiisometry.
\end{lemma}
\begin{proof}
See~\cite[Proposition~8.19 in  Chapter~I.8 on page~140]{Bridson-Haefliger(1999)}.
\end{proof}

\begin{example} \label{exa:closed_manifolds_and_pi}
Let $M = (M,g)$ be a closed connected Riemannian manifold. 
Let $\widetilde{M}$ be its universal covering. The fundamental group
$\pi = \pi_1(M)$ acts freely on $\widetilde{M}$. Equip $\widetilde{M}$
with the unique $\pi$-invariant Riemannian metric for which the projection
$\widetilde{M} \to M$ becomes a local isometry. 
The fundamental group $\pi$ is finitely generated. 
Equip it with the word metric with respect to any finite set of generators.

Then $\pi$ and $\widetilde{M}$ are quasiisometric by the
\v{S}varc-Milnor Lemma~\ref{lem:Svarc-Milnor_Lemma}. 
\end{example}

\begin{definition} \label{def:commensurable}
Two groups $G_1$ and $G_2$ are \emph{commensurable} if there are subgroups
$H_1 \subseteq G_1$ and $H_2 \subseteq G_2$ such that the indices
$[G_1:H_1]$ and $[G_2:H_2]$ are finite and $H_1$ and $H_2$ are isomorphic.
\end{definition}

\begin{lemma}\label{lem:commensurable_implies_quasiisometry}
Let $G_1$ and $G_2$ be finitely generated groups. Then:

\begin{enumerate}
\item A group homomorphism $G_1 \to G_2$ 
is a quasiisometry if and only if its kernel is finite and its image has 
finite index in $G_2$;

\item If  $G_1$ and $G_2$ are commensurable, then they are quasiisometric.

\end{enumerate}
\end{lemma}

There are quasiisometric groups that are not commensurable as the following example shows.

\begin{example} \label{exa:quasiisometry_implies_not_commnensurable}
Consider a semi-direct product $G_{\phi} = \IZ^2 \rtimes_{\phi} \IZ$
for an isomorphism $\phi \colon \IZ^2 \to \IZ^2$. For these groups
a classification up to commensurability and quasiisometry has been given
in~\cite{Bridson-Gersten(1996)} as explained next.

These groups act properly and cocompactly by isometries on precisely  one of the $3$-dimensional 
simply connected geometries
$\IR^3$, $\Nil$ or $\Sol$. (A \emph{geometry} is a complete locally homogeneous
Riemannian manifold.) If $\phi$ has finite order, then the geometry is $\IR^3$.
If $\phi$ has infinite order and the eigenvalues of the induced $\IC$-linear map
$\IC^2 \to \IC^2$ have absolute value $1$, then the geometry is $\Nil$.
If $\phi$ has infinite order and one of the eigenvalues of the induced $\IC$-linear map
$\IC^2 \to \IC^2$ has absolute value $>1$, then the geometry is $\Sol$.

These metric spaces given by the geometries $\IR^3$, $\Nil$ or $\Sol$
are mutually distinct under quasiisometry. 
By Example~\ref{exa:closed_manifolds_and_pi} two groups of the shape  
$G_{\phi}$ are quasiisometric if and only if they belong to the same geometry.

Two groups $G_{\phi}$ and $G_{\phi'}$ belonging to the same geometry $\IR^3$ or
$\Nil$ respectively contain a common subgroup of finite index and hence are commensurable.
However, suppose that $G_{\phi}$ and $G_{\phi'}$ belong to $\Sol$.
Then they are commensurable if and only if the eigenvalues $\Lambda$ and $\Lambda'$
with absolute value $> 1$
of $\phi$ and $\phi'$, respectively,
have a common power (see~\cite{Bridson-Gersten(1996)}). 
This obviously yields examples of groups
$G_{\phi}$ and $G_{\phi'}$ that belong to the geometry $\Sol$ and are quasiisometric
but are not commensurable.

The classification up to quasiisometry of finitely presented
non-poly-cyclic abelian-by-cyclic groups is presented
in~\cite[Theorem~1.1]{Farb+Mosher(2000)}.
\end{example}


\typeout{--------------------   Section 3 ---------------------------------------}

\section{Properties and invariants of groups invariant under quasiisometry}
\label{sec:Properties_and_invariants_of_groups_invariant_under_quasiisometry}

Recall that, given a property
(P) of groups, we call a group virtually-(P) if
it contains a subgroup of finite index having property (P).
In particular a group is virtually trivial if and only it is finite.
It is virtually finitely generated abelian if  and only if it 
contains a normal subgroup of finite index which is isomorphic 
to $\IZ^n$ for some integer $n \ge 0$. 

A finitely generated group $G$ is \emph{nilpotent}%
\index{group!nilpotent}
if $G$ possesses a finite \emph{lower central series}%
\index{lower central series}
$$G = G_1 \supset G_2 \supset \ldots \supset G_s = \{1\}
\hspace{15mm} G_{k+1} = [G,G_{k}].$$

A group $G$ is called \emph{amenable}
if there is a (left) $G$-invariant linear operator
$\mu\colon  l^{\infty}(G,\IR) \to \IR$
with $\mu(1) = 1$ that satisfies for all $f \in l^{\infty}(G,\IR)$
$$
\inf\{f(g)\mid g \in G\} \le
\mu(f)  \le \sup\{f(g)\mid g\in G\}.
$$
Abelian groups and more generally solvable groups are amenable.
The class of amenable groups is closed under extensions and directed unions.
A group which contains a non-abelian free group as subgroup is not amenable.
A brief survey on amenable groups
and the definition and a brief survey on accessible groups
can be found for instance in~\cite[Section~6.4.1 on page 256 ff.]{Lueck(2002)}
and~\cite[III.15 on page~52]{delaHarpe(2000)}. The book~\cite{Paterson(1988)}
is devoted to amenability. The notion of a hyperbolic
space and a hyperbolic group will be explained in
Definition~\ref{def:hyperbolic_space} and
Definition~\ref{def:hyperbolic_group}.

\begin{theorem}[Group properties invariant under quasiisometry]
\label{the:group_properties_and_quasiisometry}
The following properties of groups are \emph{geometric properties}, i.e., 
if the finitely generated group $G$ has the property, then every finitely generated
group that is quasiisometric to $G$ also has this property:

\begin{enumerate}
\item \label{the:group_properties_and_quasiisometry:finite}
Finite;
\item \label{the:group_properties_and_quasiisometry:infinite_virtually_cyclic}
Infinite virtually cyclic;
\item \label{the:group_properties_and_quasiisometry:finitely_presented}
Finitely presented;
\item \label{the:group_properties_and_quasiisometry:virtually_abelian}
Virtually abelian;
\item \label{the:group_properties_and_quasiisometry:virtually_nilpotent}
Virtually nilpotent;
\item \label{the:group_properties_and_quasiisometry:virtually_free}
Virtually free;
\item \label{the:group_properties_and_quasiisometry:amenable}
Amenable;
\item \label{the:group_properties_and_quasiisometry:hyperbolic}
Hyperbolic;
\item \label{the:group_properties_and_quasiisometry:accessible}
Accessible;
\item \label{the:group_properties_and_quasiisometry:finite_n_skeleton}
The existence of a model for the classifying space $BG$ 
with finite $n$-skeleton for given $n \ge 2$;
\item \label{the:group_properties_and_quasiisometry:finite_type}
The existence of a model for $BG$ of finite type, i.e., all skeletons are finite.
\end{enumerate}
\end{theorem}
\begin{proof}%
\ref{the:group_properties_and_quasiisometry:finite} 
Having bounded diameter is a quasiisometry invariant of metric spaces.
\\[1mm]%
\ref{the:group_properties_and_quasiisometry:infinite_virtually_cyclic}
This follows from 
Theorem~\ref{the:end_of_groups}~\ref{the:end_of_groups:2} and
Theorem~\ref{the:quasiisometry_invariants}~%
\ref{the:quasiisometry_invariants:number_of_ends}.
\\[1mm]%
\ref{the:group_properties_and_quasiisometry:finitely_presented}
See~\cite[Proposition~4 In Chapter V.A on page~119]{delaHarpe(2000)}.
\\[1mm]%
\ref{the:group_properties_and_quasiisometry:virtually_abelian}
See~\cite[Chapter~I]{Ghys-Harpe(1990)}.
\\[1mm]%
\ref{the:group_properties_and_quasiisometry:virtually_nilpotent}
This follows from 
Theorem~\ref{the:quasiisometry_invariants}~\ref{the:quasiisometry_invariants:growth}
and Theorem~\ref{the:virtually_nilpotent_and_growth}.
\\[1mm]%
\ref{the:group_properties_and_quasiisometry:virtually_free}
See~\cite[Theorem~19 in Chapter~I]{Ghys-Harpe(1990)} and
Theorem~\ref{the:quasiisometry_invariants}%
\ref{the:quasiisometry_invariants:number_of_ends}.
\\[1mm]%
\ref{the:group_properties_and_quasiisometry:amenable}
This follows from~\cite{Foelner(1955)}.
See also~\cite[Chapter~6]{Gromov(1981)}.
\\[1mm]%
\ref{the:group_properties_and_quasiisometry:hyperbolic}
Quasiisometric groups have 
quasiisometric Cayley graphs and it is not difficult to see that
the property being hyperbolic is a quasiisometry invariant of geodesic spaces.
\\[1mm]%
\ref{the:group_properties_and_quasiisometry:accessible}
See~\cite{Thomasson+Woess(1993)}.
\\[1mm]%
\ref{the:group_properties_and_quasiisometry:finite_n_skeleton}
This follows from~\cite[Item~1.$C_2^{\prime}$ on page~25]{Gromov(1993)}.
See also~\cite{Alonso(1994)}.
\\[1mm]%
\ref{the:group_properties_and_quasiisometry:finite_type} 
This follows from assertion~\ref{the:group_properties_and_quasiisometry:finite_n_skeleton}.
\end{proof}

If $S$ is a finite set of generators for the group $G$,
let $b_S(n)$ be the number of elements in $G$ which can be written
as a word in $n$ letters of $S \cup S^{-1} \cup \{1\}$, i.e., the number of elements
in the closed ball of radius $n$ around $1$ with respect to $d_S$. 

The following definition is indeed independent  of the choice of
the finite set $S$ of generators.
\begin{definition}[Growth]\label{def:polynomial_growth}
The group $G$ has
\emph{polynomial growth of degree not greater than $d$}
if there is $C$ with $b_S(n) \le C n^d$ for all $n \ge 1$.

We say that $G$ has \emph{polynomial growth}
if it has polynomial growth of degree not greater than $d$ for some $d > 0$.

It has \emph{exponential growth} if there exist $C > 0$ and $\alpha > 0$ 
such that for $n \ge 1$ we have
$$b_S(n) \ge C \cdot \alpha^n.$$
It has \emph{subexponential growth} if it has neither 
polynomial growth nor exponential growth.
\end{definition}

The free abelian group $\IZ^n$ of rank $n$ has polynomial growth rate
of precisely degree $n$. A finitely generated non-abelian free  group has
exponential growth rate.

Recall that the \emph{Hirsch rank} of a solvable group $G$ is defined to be
$$h(G) = \sum_{i \ge 0} \dim_{\IQ}(G_{i+1}/G_i \otimes_{\IZ} \IQ),$$
where $G_i$ is the $i$-th term in the derived series of $G$.

A metric is called \emph{proper} if every closed ball is compact.
Let $X$ be a proper geodesic space. A \emph{proper ray} is a map 
$[0,\infty) \to X$ such that the preimage of a compact set is compact again.
Two proper rays $c_0,c_1 \colon [0,\infty) \to X$
\emph{converge to the same end} if for every compact subset $C \subset X$ there is 
$R > 0$ such that $c_0([R,\infty))$  and $c_1([R,\infty))$ lie in the same component
of $X\setminus C$. This defines an equivalence relation on the set of proper rays.
The set of equivalence classes is the \emph{set of ends}. 

The \emph{number of ends} of $X$ is the
cardinality of this set. It is a quasiisometry invariant
(see~\cite[Proposition~8.29 on page~128]{Bridson-Haefliger(1999)}).
Hence the following definition makes sense.

\begin{definition}[Number of ends]
\label{def:number_of_ends}
The \emph{number of ends} of a finitely generated group $G$ is defined to be
the number of ends of the Cayley graph $\Cay_S(G)$ for any choice 
of a finite set $S$  of generators. 
\end{definition}

\begin{theorem}[Ends of groups]\
\label{the:end_of_groups}
\begin{enumerate}
\item \label{the:end_of_groups:possible_numbers}
A finitely generated group has $0$, $1$, $2$ or infinitely many ends;

\item \label{the:end_of_groups:0}
It has $0$ ends precisely if it is finite;

\item \label{the:end_of_groups:2}
It has two ends precisely if it is infinite and virtually cyclic.

\item \label{the:end_of_groups:infinity}
It has infinitely many ends if and only if $G$ can be expressed as an amalgamated
product $A \ast_C B$ or as an $\HNN$-extension $A \ast_C$ with finite $C$ and
$|A/C| \ge 3$ and $|B/C| \ge 2$. 
\end{enumerate}
\end{theorem}
\begin{proof}
See~\cite[Theorem~8.32 in Chapter~I.8 on page~146]{Bridson-Haefliger(1999)}.
\end{proof}

\begin{theorem}[Invariants under quasiisometry]
\label{the:quasiisometry_invariants}
Let $G_1$ and $G_2$ be two finitely generated groups which are quasiisometric. Then:
\begin{enumerate}
\item \label{the:quasiisometry_invariants:number_of_ends}
They have the same number of ends;

\item \label{the:quasiisometry_invariants:cohomological_dimension}

Let $R$ be a commutative ring. Then we get
$$\cd_R(G_1) = \cd_R(G_2)$$
if one of the following assumptions is satisfied:
\begin{enumerate}
\item \label{the:quasiisometry_invariants:cohomological_dimension:(1)} 
The cohomological dimensions $\cd_R(G_1)$ and $\cd_R(G_2)$ are both finite;

\item \label{the:quasiisometry_invariants:cohomological_dimension:(2)} 
There exist finite models for $BG_1$ and $BG_2$;

\item \label{the:quasiisometry_invariants:cohomological_dimension:(3)}
One of the groups $G_1$ and $G_2$ is amenable and $\IQ \subseteq R$;

\end{enumerate}

\item \label{the:quasiisometry_invariants:Hirsch_length}
If they are solvable, then they have the same Hirsch length;

\item \label{the:quasiisometry_invariants:growth}
Suppose that $G_1$ has polynomial growth of degree not greater than $d$,
intermediate growth,or  exponential growth, respectively. Then the same is true
for $G_2$;

\item \label{the:quasiisometry_invariants:real_cohomology_ring}
Let $G_1$ and $G_2$ be nilpotent.
Then their real cohomology rings $H^*(G_1;\IR)$ and  $H^*(G_2;\IR)$
are isomorphic as graded rings. In particular the Betti numbers
of $G_1$ and $G_2$ agree.

\end{enumerate}
\end{theorem}
\begin{proof}%
\ref{the:quasiisometry_invariants:number_of_ends}
See~\cite[Corollary~2.3]{Brick(1993)} or~\cite[Corollary~1]{Gersten(1993)}.
\\[1mm]%
\ref{the:quasiisometry_invariants:cohomological_dimension}
See~\cite[Theorem~1.2]{Sauer(2006)}. The case $R = \IZ$ under 
condition~\eqref{the:quasiisometry_invariants:cohomological_dimension:(2)}
has already been treated in~\cite[Corollary~2]{Gersten(1993)}.
\\[1mm]%
\ref{the:quasiisometry_invariants:Hirsch_length}
This follows from assertion~\ref{the:quasiisometry_invariants:cohomological_dimension}
since $\cd_{\IQ}(G)$ is the Hirsch rank for a virtually poly-cyclic group $G$
(see~\cite[Corollary~1.3]{Sauer(2006)}). 
\\[1mm]%
\ref{the:quasiisometry_invariants:growth}
See~\cite[Proposition~27 in VI.B on page 170]{delaHarpe(2000)}.
\\[1mm]%
\ref{the:quasiisometry_invariants:real_cohomology_ring}
See~\cite[Theorem~1.5]{Sauer(2006)}. The statement about the Betti
numbers was already proved by Shalom~\cite[Theorem~1.2]{Shalom(2004)}.

\end{proof}

We mention that there is an extension of the notion of quasiisometry
to groups which are not necessarily finitely generated and that some of the
the results of Theorem~\ref{the:quasiisometry_invariants} are still true in this more 
general setting (see~\cite{Sauer(2006)}, \cite{Shalom(2004)}).

\begin{conjecture}[Folk]\label{con:folk_conjecture}
Let $G_1$ and $G_2$ be two finitely generated torsionfree nilpotent groups. 
Let $L_1$ and $L_2$ be the simply connected nilpotent Lie groups given by their Mal'cev
completion. (These are uniquely determined by the fact that $G_i$ is cocompactly embedded
in $L_i$.)

If $G_1$ and $G_2$ are quasi-isometric, then $L_1$ and $L_2$
are isomorphic as Lie groups.
\end{conjecture}

\begin{remark} \label{rem:Pansu}
Evidence for Conjecture~\ref{con:folk_conjecture} comes from the following facts.
The graded Lie algebra associated to the Mal'cev completion of a finitely generated
torsionfree nilpotent group $G$ is a quasiisometry invariant of $G$
by a result of Pansu~\cite{Pansu(1989c)}.
The result of Sauer mentioned in 
Theorem~\ref{the:quasiisometry_invariants}~%
\ref{the:quasiisometry_invariants:real_cohomology_ring}
follows from Conjecture~\ref{con:folk_conjecture} since the 
cohomology algebras of the Lie algebra of the Mal'cev  completion and the cohomology algebra
of $G$ itself are isomorphic (see~\cite[Theorem~1]{Nomizu(1954)}).
\end{remark}

The following celebrated theorem due to Gromov~\cite{Gromov(1981c)} is one of the milestones
in geometric group theory. A new proof can be found in~\cite{Kleiner(2007)}.

\begin{theorem}[Virtually nilpotent groups and growth]
\label{the:virtually_nilpotent_and_growth}
A finitely generated group is virtually nilpotent if and only if it has polynomial growth.
\end{theorem}

\begin{remark}[Virtually solvable groups]
\label{rem:virtually_solvable_groups}
This raises the question whether
solvability is a geometric property. However, there exists a finitely
generated solvable group which is 
quasiisometric to a finitely generated group which is not virtually 
solvable (see~\cite{Dioubina(1999)}).
This counterexample is not finitely presented.
It is still not known whether two finitely presented quasiisometric groups both
have to be virtually solvable if one of them is. 
\end{remark}

\begin{remark}[Free products] \label{free_products}
Let $G_1$, $G_1'$, $G_2$ and $G_2'$ be finitely generated groups.
Suppose that $G_i$ and $G_i'$ are quasiisometric for $i = 1,2$.
Assume that none of the groups $G_1$, $G_1'$, $G_2$ and $G_2'$ is trivial or $\IZ/2$.
Then the free products $G_1 \ast G_2$ and $G_1' \ast G_2'$ are
quasiisometric. (They are actually Lipschitz equivalent).
See~\cite[46.~(ii) in IV.B on page 105]{delaHarpe(2000)}
and~\cite[Theorem~0.1]{Papasoglu+Whyte(2002)}.
The corresponding statement is false if one replaces
quasiisometric by commensurable (see~\cite[46.~(iii) in IV.B on page 106]{delaHarpe(2000)}).
\end{remark}

\begin{remark}[Property (T)]
\label{rem:property_(T)}
Kazhdan's Property (T) is not a quasiisometry invariant.
(This is due to Furman and Monod and stated in~\cite[page~173]{Gaboriau(2002b)}).
\end{remark}

\begin{remark}[The sign of the Euler characteristic]
\label{rem:sign_of_Euler:_characteristic}
The sign of the Euler characteristic of a group with a finite model
for $BG$ is not a quasiisometry invariant. 
See~\cite[46.~(iii) in IV.B on page 105]{delaHarpe(2000)}.
\end{remark}

\begin{remark}[Minimal dimension of $EG$ and $\eub{G}$]
\label{rem:minimal_dimension_of_eub(G)}
We have already mentioned in Theorem~\ref{the:quasiisometry_invariants}~%
\ref{the:quasiisometry_invariants:cohomological_dimension}
that the cohomological dimension $\cd_{\IZ}(G)$ is a quasiisometry
invariant under the assumption that there exists a $G$-$CW$-model for $EG$
which is finite or, equivalently, cocompact.

There always exists
a $\max\{3,\cd_{\IZ}(G)\}$-dimensional model for $BG$
(see~\cite[Theorem~7.1 in Chapter~VII on page 295]{Brown(1982)}).
Notice that the existence of a $d$-dimensional $CW$-model for $BG$ is
equivalent to the existence of a $d$-dimensional $G$-$CW$-model for $EG$
since $EG$ is the universal covering of $BG$.
Hence $\cd_{\IZ}(G)$  is equal to the minimal dimension of a 
model for $EG$ if $\cd_{\IZ}(G) \ge 3$.

If $H \subset G$ is a subgroup of finite index of the torsionfree group $G$
and there is a finite dimensional
model for $EG$, then the cohomological dimensions of $G$ and $H$ agree
by a result of Serre (see~\cite[Theorem~3.2 in Chapter~VIII.3  on page~190]{Brown(1982)}, 
\cite{Serre(1971)}) and hence also the minimal dimension for $EH$ and $EG$
agree if the cohomological dimension of $G$ is greater or equal to $3$.

The corresponding statement is false if one replaces $EG$ by the universal
space $\eub{G}$ for proper group $G$-actions 
(see Definition~\ref{def:classifying_space_for_proper_actions}). 
Namely, there exists a group
$G$ with  a torsionfree subgroup $H$ of finite index such that there exists
a $d$-dimensional model for $EH = \eub{H}$ but no
$d$-dimensional model for $\eub{G}$ (see~\cite[Theorem~6]{Leary-Nucinkis(2003)}).

Hence the minimal dimension of a model for $\eub{G}$ is not at all a quasiisometry invariant
in general. 
\end{remark}

\begin{remark}[$L^2$-invariants]
\label{rem:L2-invariants}
If the finitely generated groups $G_1$ and $G_2$
are quasiisometric and there exist finite models for $BG_1$ and $BG_2$ then 
$$b_p^{(2)}(G_1) = 0 \Leftrightarrow b_p^{(2)}(G_2)= 0$$ 
holds (see \cite[page 224]{Gromov(1993)}, \cite{Pansu(1995)}).
But it is general not true that in the situation above there exists
a constant $C > 0$ such that $b_p^{(2)}(G_1) = C \cdot  b_p^{(2)}(G_2)$ holds
for all $p \ge 0$ (see~\cite[page~313]{Lueck(2002)}, \cite{Papasoglu(1995trees)}).

It is unknown whether the vanishing of the $L^2$-torsion of appropriate
groups or the Novikov-Shubin invariants of appropriate groups are quasiisometry invariants.
(see~\cite[Question~7.35 and Question~7.36 on page~313]{Lueck(2002)}).
Partial results in this direction have been obtained 
in~\cite[Theorem~1.6]{Sauer(2006)} and in~\cite{Wegner(2001)} for amenable respectively
elementary amenable groups.
\end{remark}

\begin{remark}[Asymptotic cone] \label{rem:asymptotic_cone}
The notion of an asymptotic cone using ultralimits
was introduced by Van den Dries and Wilkie~\cite{Dries+Wilkie(1984)}.
It assigns to a metric  space a new space after the choice
of a non-principal ultrafilter on the set of natural numbers, 
a scaling sequence and a sequence of observation points.
The asymptotic cone does in general depend on these extra choices.
Roughly speaking, an asymptotic cone of a metric space is what 
one sees when one looks at the space from infinitely far away. 

Applied to the Cayley graph of a finitely generated group
an asymptotic cone yields a complete geodesic homogeneous metric space,
which captures the coarse properties. It depends on the ultrafilter and 
the scaling sequence but not on the sequence of observation points. 
A quasiisometry induces a bi-Lipschitz homeomorphism between the asymptotic cones
(for the same ultrafilters and scaling constants). 
So as in the case of the boundary of a hyperbolic group 
(see Section~\ref{sec:The_boundary_of_a-hyperbolic_spaces}) 
we can assign to a group a metric
space such that a quasiisometry induces a ``nice'' map between the associated structures.

Further information and a discussion of some applications to 
quasiisometry can be found for instance
in~\cite[Chapter~I.5 on page 77 ff.]{Bridson-Haefliger(1999)} 
and~\cite{Drutu(2002)}. 
Asymptotic cones play a significant role
in the proof of certain rigidity results, for instance in  the proof
of the rigidity of quasiisometries for symmetric spaces and Euclidean buildings
due to Kleiner-Leeb~\cite{Kleiner-Leeb(1997)} or in the 
proof the rigidity under quasiisometry of the mapping class group
(see Theorem~\ref{the:quasiisometric_rigidity_mapping_class_group}) 
due to Behrstock-Kleiner-Minsky-Mosher~\cite{Behrstock-Kleiner-Minsky-Mosher(2008)}
and Hamenst\"adt~\cite{Hamenstaedt(2007)}. Asymptotic cones and quasiisometry classes
of fundamental groups of $3$-manifolds are investigated 
in~\cite{Kapovich-Leeb(1995)}.
\end{remark}

\begin{remark}[Group splittings] \label{rem:Group_splitting}
A lot of activity in geometric group theory has been focused on extending
the Jaco-Johannson-Shalen decomposition for $3$-manifolds to finitely presented groups
(see for instance~\cite{Bowditch(1998)}, \cite{Dunwoody+Sageev(1999)},
\cite{Rips-Sela(1997)}, \cite{Scott+Swarup(2003)}, \cite{Sela(1997)}).
Its quasiisometry invariance has been proved in~\cite{Papasoglu(2005)}.
\end{remark}

Further information about quasiisometry invariants can be found for instance
in~\cite{Block+Weinberger(1997)}, \cite{Bridson-Haefliger(1999)}, \cite{delaHarpe(2000)},
\cite{Ghys+delaHarpe(1991)}, \cite{Gromov(1993)}.


\typeout{--------------------   Section 4 ---------------------------------------}

\section{Rigidity}
\label{sec:Rigidity}

An explanation of the following two theorems and a 
list of papers that have made significant contributions to their proof
can be found in~\cite{Farb(1997)}. It includes
Eskin~\cite{Eskin(1998)},
Eskin-Farb~\cite{Eskin+Farb(1998)},
Farb-Schwarz~\cite{Farb+Schwartz(1996)}, Kleiner-Leeb~\cite
{Kleiner-Leeb(1997)}, Pansu~\cite {Pansu(1989c)},
and Schwartz~\cite{Schwartz(1995)} and \cite{Schwartz(1996)}.

In the sequel \emph{semisimple Lie group}
means non-compact, connected semisimple Lie group with finite center.
\emph{Lattice} means a discrete subgroup of finite covolume.
A lattice is called \emph{uniform} if it is cocompact.

\begin{theorem}[Rigidity of the class of lattices] 
\label{the:rigidity_of_the_class_of_lattices}
Let $\Gamma$ be a finitely generated group. If $\Gamma$ is quasiisometric to an irreducible
lattice in a semisimple Lie group $G$, then $\Gamma$ is almost a lattice in $G$, i.e.,
there is a lattice $\Lambda'$ in $G$ and a finite group $F$ such that there exists
an exact sequence
$$1 \to F \to \Gamma \to \Lambda' \to 1.$$
\end{theorem}

\begin{theorem}[Classification among lattices]
\label{the:classification_among_lattices}
The quasiisometry classes of irreducible lattices in semisimple Lie groups are precisely:

\begin{enumerate}

\item \label{the:classification_among_lattices:uniform}
One quasiisometry class for each semisimple Lie group, consisting of the
uniform lattices in $G$;

\item \label{the:classification_among_lattices:non-uniform}
One quasiisometry class for each commensurability class of irreducible non-uniform
lattices, except in  $G = \mathrm{SL}_2(\IR)$, where there is precisely 
one quasiisometry class of 
non-uniform lattices.

\end{enumerate}
\end{theorem}

The following result is the main result of~\cite{Farb+Mosher(1999)}.
Recall that for $n \ge 2$ the solvable \emph{Baumslag-Solitar-group} is defined by
$$\mathrm{}BS(1,n) = \langle a,b \mid bab^{-1} = a^n\rangle.$$

\begin{theorem}[Rigidity of Baumslag-Solitar groups]
\label{the:solvable_Baumslag-Solitar-group}
Let $G$ be a finitely generated group. Suppose that $G$ is quasiisometric to
$\mathrm{BS}(1,n)$ for some $n \ge 2$. Then there is an exact sequence
$$1 \to F \to G \to \Gamma \to 1,$$
where $F$ is finite and $\Gamma$ is commensurable to $BS(1,n)$.
\end{theorem}

\begin{remark}[Abelian by cyclic groups] \label{rem:abelian_by_cyclic}
The quasiisometry rigidity of finitely presented abelian-by-cyclic groups is investigated
in~\cite[Theorem~1.2]{Farb+Mosher(2000)}.
\end{remark} 

The following result is due 
Behrstock-Kleiner-Minsky-Mosher~\cite{Behrstock-Kleiner-Minsky-Mosher(2008)}
and Hamen\-st\"adt~\cite{Hamenstaedt(2007)}.

\begin{theorem}[Rigidity of mapping class groups] 
\label{the:quasiisometric_rigidity_mapping_class_group}
Let $S$ be an oriented closed surface. 
Let  $M(S)$ be the associated mapping class group. 
Let $G$ be a finitely generated group that is quasiisometric to
$M(S)$.
Let $\cent(M(S))$ be the center of $M(S)$ which is a finite group.

Then there is a finite index
subgroup $G'$ in $G$ and a homomorphism $G'  \to M(S)/\cent(M(S))$ with finite kernel and
finite index image.
\end{theorem}


\typeout{--------------------   Section 5 ---------------------------------------}

\section{Hyperbolic  spaces and CAT($\kappa$)-spaces}
\label{sec:Hyperbolic_spaces_and:cat(kappa)-spaces}

Recall that we have introduced the notion of a geodesic space in 
Definition~\ref{def:geodesic_space}.

\begin{example}[Geodesic spaces] \label{exa:geodesic_spaces}
A  complete Riemannian manifold 
inherits the structure of a geodesic metric space from the Riemannian metric
by defining the distance of two points to be the infimum over the length of any curve joining 
them.

A graph inherits the structure of a metric space by defining the distance of two points
to be the infimum over the length of any piecewise linear path joining them,
where each edge is isometrically identified with the unit interval $[0,1]$.
A graph is connected if and only if it is a geodesic space with respect to this metric.
\end{example}

A \emph{geodesic triangle} in a geodesic space $X$ is a configuration of
three points $x_1$, $x_2$ and $x_3$  in $X$ together with a choice of 
three geodesics $g_1$, $g_2$ and $g_3$ such that $g_1$ joins $x_2$ to $x_3$,
$g_2$ joins $x_1$ to $x_3$ and $g_3$ joins $x_1$ to $x_2$. For $\delta > 0$
a geodesic triangle is called \emph{$\delta$-thin} if each edge is contained in
the closed $\delta$-neighborhood of the union of the other two edges.

\begin{definition}[Hyperbolic space] \label{def:hyperbolic_space}
Consider $\delta \ge 0$. A \emph{$\delta$-hyperbolic space} is a  geodesic space
whose geodesic triangles are all $\delta$-thin.

A geodesic space is called \emph{hyperbolic} it is $\delta$-hyperbolic for some
$\delta > 0$.
\end{definition}

\begin{remark}[Equivalent definitions of hyperbolic space]
\label{rem:equivalent_definitions_of_hyperbolic_space}
There are many equivalent definitions of hyperbolic spaces, 
which are useful and can be found under the key words
``fine triangles'', ``minsize'', ``insize'', ``Gromov's inner
product and the $4$-point-condition'', ``geodesic divergence'' and ``linear isoperimetric
inequality''.
(see for instance~\cite{Bowditch(1991)}, \cite{Bridson-Haefliger(1999)},
\cite{Ghys-Harpe(1990)}, \cite{Gromov(1987)}).
\end{remark}

\begin{remark}[Examples and non-examples for hyperbolic spaces] 
\label{rem:examples_ofhyperbolic_spaces}\
Every geodesic space with bounded diameter is hyperbolic.
Every complete Riemannian manifold whose sectional curvature is bounded from above by a
negative constant is a hyperbolic space. In particular the hyperbolic $n$-space
$\IH^n$ and every closed Riemannian manifold with negative sectional curvature
are hyperbolic spaces. The Euclidean space $\IR^n$ is not hyperbolic.
A tree is $\delta$-hyperbolic for every $\delta \ge 0$.
\end{remark}

For $\kappa \le 0$ let $M_{\kappa}$ be the up to isometry unique simply connected complete
Riemannian manifold whose sectional curvature is constant with value $\kappa$. Consider
a metric space $X$. For every geodesic triangle $\Delta$ with edges $x_1$, $x_2$ and $x_3$
in $X$ there exists a geodesic
triangle $\overline{\Delta}$ in $M_{\kappa}$ with edges $\overline{x_1}$, $\overline{x_2}$,
and $\overline{x_3}$ which is a geodesic triangle and satisfies
$d_X(x_i,x_j) = d_{M_{\kappa}}(\overline{x_i},\overline{x_j})$ for $i,j \in \{1,2,3\}$.
We call such a triangle $\overline{\Delta}$ a \emph{comparison triangle}.
It is unique up to isometry. For every point $x$ in  $\Delta$ there is unique
comparison point $\overline{x}$ determined by the property that $\overline{x}$  lies
on the edge from $\overline{x_i}$ to $\overline{x_j}$ if $x$ lies on the edge
joining $x_i$ and $x_j$ and the distance of $\overline{x}$ and
$\overline{x_i}$ agrees with the distance of $x$ and $x_i$.

\begin{definition}[$\CAT(\kappa)$-space]
\label{def:cat(kappa)_space}
Let $X$ be a geodesic space and let $\kappa \le 0$. Then $X$ satisfies
the \emph{$\CAT(\kappa)$-condition} if for every geodesic triangle $\Delta$ and points
$x,y \in \Delta$ and any comparison triangle $\overline{\Delta}$ in $M_{\kappa}$
and comparison points $\overline{x}$ and $\overline{y}$ we have
$$ d_X(x,y) \le d_{M_{\kappa}}(\overline{x},\overline{y}).$$

A \emph{$\CAT(\kappa)$-space} is a geodesic space which  
satisfies the $\CAT(\kappa)$-condition.

A geodesic space  is \emph{of curvature $\le \kappa$} for some $\kappa \le 0$
if it satisfies the $\CAT(\kappa)$-condition locally. It is called 
\emph{negatively curved} or \emph{non-positively curved} respectively
if it is of curvature $\le \kappa$ for some $\kappa < 0$ or $\kappa \le 0$ respectively.
\end{definition}

A space $Y$ is called \emph{aspherical} if it is path connected and 
$\pi_n(Y,y)$ vanishes for one (and hence all) $y \in Y$. 
Provided that $Y$ is a $CW$-complex, $Y$ is aspherical
if and only if it is connected and its universal covering is contractible.

\begin{theorem}[$\CAT(\kappa)$-spaces] \label{the:CAT(kappa)-spaces}
Fix $\kappa \le 0$. Then:
\begin{enumerate}

\item \label{the:CAT(kappa)-spaces:Riemann}
A simply connected Riemannian manifold has sectional curvature $\le \kappa$ 
if and only if it is a $\CAT(\kappa)$-space with respect to the metric induced 
by the Riemannian metric;

\item \label{the:CAT(kappa)-spaces:contractible}
A $\CAT(\kappa)$-space is contractible;

\item \label{the:CAT(kappa)-spaces:simply_connected}
A simply connected complete geodesic space of curvature $\le \kappa$ is a $\CAT(\kappa)$-space;

\item \label{the:CAT(kappa)-spaces:aspherical}
A complete geodesic space of curvature $\le \kappa$ has a $\CAT(\kappa)$-space 
as universal covering and is aspherical;

\item \label{the:CAT(kappa)-spaces:kappa_le_kappaprime}
Consider $\kappa \le \kappa' \le 0$. If $X$ is a $\CAT(\kappa)$-space 
of curvature $\le \kappa$, then $X$ is a $\CAT(\kappa')$-space 
of curvature $\le \kappa'$;

\item \label{the:CAT(kappa)-spaces:hyperbolic}
A proper $\CAT(0)$-space is hyperbolic if and only if it contains
no subspace isometric to $\IR^2$;

\item \label{the:CAT(kappa)-spaces:kappa_negative}
For $\kappa < 0$ a $\CAT(\kappa)$-space is hyperbolic;

\item \label{the:CAT(kappa)-spaces:tree}
A tree is a $\CAT(\kappa)$-space for all $\kappa \le 0$.

\end{enumerate}
\end{theorem}
\begin{proof}%
\ref{the:CAT(kappa)-spaces:Riemann}
See~\cite[Corollary~1A.6 in Chapter~II.1 on page~173]{Bridson-Haefliger(1999)}.
\\[1mm]%
\ref{the:CAT(kappa)-spaces:contractible}
See~\cite[Corollary~1.5 in Chapter~II.1 on page~161]{Bridson-Haefliger(1999)}.
\\[1mm]%
\ref{the:CAT(kappa)-spaces:simply_connected}
See~\cite[Theorem~4.1~(2) in Chapter~II.4 on page~194]{Bridson-Haefliger(1999)}.
\\[1mm]%
\ref{the:CAT(kappa)-spaces:aspherical}
This follows from
assertions~\ref{the:CAT(kappa)-spaces:contractible} and~\ref{the:CAT(kappa)-spaces:simply_connected}.
\\[1mm]%
\ref{the:CAT(kappa)-spaces:kappa_le_kappaprime}
See~\cite[Theorem~1.12 in Chapter~II.1 on page~165]{Bridson-Haefliger(1999)}.
\\[1mm]%
\ref{the:CAT(kappa)-spaces:hyperbolic}
See~\cite[Theorem~1.5 in Chapter~III.H on page~400]{Bridson-Haefliger(1999)}.
\\[1mm]%
\ref{the:CAT(kappa)-spaces:kappa_negative}
See~\cite[Proposition~1.2 in Chapter~III.H on page~399]{Bridson-Haefliger(1999)}.
\\[1mm]%
\ref{the:CAT(kappa)-spaces:tree}
See~\cite[Example~1.15~(4)  in Chapter~II.1 on page~167]{Bridson-Haefliger(1999)}.
\end{proof}

\begin{remark} \label{rem:hyperbolic_versus_CAT}
The condition of being
hyperbolic is a condition in the large. For instance, a local change of the metric 
on a compact subset does not destroy this property. 
This is not true for the condition being $\CAT(\kappa)$. For example,
any compact metric space is hyperbolic, whereas it is
not $\CAT(\kappa)$ for some $\kappa \le 0$  in general.  

In general it makes a significant difference whether a space is negatively curved
or non-positively curved.

There is no version of the $\CAT(0)$-condition known that
is like the condition hyperbolic defined in the large.
\end{remark}


\typeout{--------------------   Section 6 ---------------------------------------}

\section{The boundary of a hyperbolic  space}
\label{sec:The_boundary_of_a-hyperbolic_spaces}

Let $X$ be a hyperbolic space. 
A \emph{geodesic ray} is a geodesic $c \colon [0,\infty) \to X$ with 
$[0,\infty)$ as source.
We call two geodesic rays $c,c' \colon [0,\infty) \to X$ \emph{asymptotic} if there exists
$C \ge 0$ such that $d_X(c(t),c'(t)) \le C$ holds for all $t \in [0,\infty)$.

\begin{definition}[Boundary of a hyperbolic space]
\label{def:boundary}
Let $\partial X$ be the set of equivalence classes of geodesic rays.
Put
$$\overline{X} = X \amalg \partial X.$$
\end{definition}

The description of the topology on $\overline{X}$ and the proof of the following
two results can be found in~\cite[Chapter~III.H on pages~429-430 and
Exercise~3.18~(4) in Chapter~III.H on page~433]{Bridson-Haefliger(1999)}.

\begin{lemma}
\label{lem:topology_on_overline_X}
There is a topology on $\overline{X}$ such that $\overline{X}$ is compact and metrizable,
the subspace topology of $X \subseteq \overline{X}$ agrees with the topology coming from
the metric, $X \subseteq \overline{X}$ is open and dense, and 
$\partial X \subseteq \overline{X}$ is closed.
\end{lemma}

\begin{lemma} \label{lem:quasiisometry_and_boundary}
Let $X$ and $Y$ be hyperbolic spaces. Let $f \colon X \to Y$ be a quasiisometry.
It induces a map 
$$\overline{f} \colon \overline{X} \to \overline{Y},$$
which restricts on the boundary to a homeomorphism
$$\partial f \colon \partial X \xrightarrow{\cong} \partial Y.$$
In particular, the boundary is a quasiisometry invariant of a hyperbolic space.
\end{lemma}

\begin{remark}[Mostow rigidity] \label{rem:Mostow_rigidity}
Let $f \colon M \to N$ be a homotopy equivalence of hyperbolic closed manifolds
of dimension $n \ge 3$. \emph{Mostow rigidity} says that $f$ is homotopic to an isometric
diffeomorphism. Lemma~\ref{lem:quasiisometry_and_boundary} plays a role in its proof
as we briefly explain next. More details can be found for instance
in~\cite{Benedetti-Petronio(1992)}.

Notice that the universal coverings $\widetilde{M}$ and $\widetilde{N}$
are isometrically diffeomorphic to the $n$-dimensional hyperbolic space
$\IH^n$. The boundary of $\IH^n$ can be identified with 
$S^{n-1}$ and $\overline{\IH^n}$ with $D^n$.
Since  $M$ and $N$ are compact,
the map $\widetilde{f} \colon \widetilde{M} \to \widetilde{N}$ is a quasiisometry. 
Hence it induces a homeomorphism $\partial  \widetilde{f} \colon \partial \widetilde{M}
\xrightarrow{\cong} \partial \widetilde{N}$.  
Next one shows that the volume of a closed hyperbolic manifold is a homotopy invariant, for instance
using the notion of the simplicial volume due to Gromov and Thurston.
This is used to prove that an ideal simplex in $\overline{\IH^n}$ 
with vertices $x_0$, $x_1$, \ldots, $x_n$ on $\partial{\IH^n}$
has the same volume as the ideal triangle with vertices $\partial \widetilde{f}(x_0)$, 
$\partial \widetilde{f}(x_1)$, \ldots, $\partial \widetilde{f}(x_{n})$. This implies
that there is an isometric diffeomorphism $\widetilde{g} \colon \widetilde{M} \to \widetilde{N}$
with $\partial \widetilde{g} = \partial \widetilde{f}$ such that 
$\widetilde{g}$ is compatible with the actions of the fundamental groups and passes to an isometric
diffeomorphism $g \colon M \to N$ which induces on the fundamental groups the same map as $f$ and hence
is homotopic to $f$. 

In the last step the condition $n \ge 3$ enters. Indeed, 
Mostow rigidity does not hold in dimension $n = 2$.
\end{remark}


\typeout{--------------------   Section 7 ---------------------------------------}

\section{Hyperbolic groups}
\label{sec:Hyperbolic_groups}

\begin{definition}[Hyperbolic group] \label{def:hyperbolic_group}
A finitely generated group is called \emph{hyperbolic} if its Cayley graph
is hyperbolic.
\end{definition}

Recall that the quasiisometry type of the Cayley graph of a
finitely generated group $G$ depends only on
$G$ as a group but not on the choice of a finite set of generators and
the notion hyperbolic is a quasiisometry invariant for geodesic
spaces. Hence the definition above makes sense, i.e., being hyperbolic is
a property of the finitely generated group $G$ itself and does not depend on the choice
of a finite set of generators. 

Let $G$ be hyperbolic. Its \emph{boundary}
$\partial G$  is the boundary of the Cayley graph.
This is well-defined up to homeomorphism, i.e., independent of a choice of a finite set of generators
because of Lemma~\ref{lem:quasiisometry_and_boundary}.

The notion of the \emph{classifying space for proper $G$-actions}
$\eub{G}$  will be explained in Definition~\ref{def:classifying_space_for_proper_actions}.

A \emph{Dehn  presentation} of a  group $G$ 
with a finite set of generators $S$ is
a finite list of words $u_1$, $v_1$, \ldots , $u_n$, $v_n$ such that
$u_1 = v_1$, \ldots , $u_n = v_n$ holds in $G$, and $d_S(e,v_i) \le
d_S(e,u_i)$ is true for $i = 1,2, \ldots , n$ and any word
represents the identity element $e$ only if it contains one of the
words $u_i$ as a subword. Now there is an obvious algorithm to decide whether a word
$w$ represents the unit element $e$ in $G$: Look  whether it contains
one of the words $u_i$. If the answer is no, the process ends, if the
answer is yes, replace $u_i$ by $v_i$. By induction over $d_S(e,w)$ one sees
that this process stops after at most $d_S(e,w)$  steps.
The word $w$ represents $e$ if and only if  the process ends
with the trivial word.

A survey article about  \emph{Poincar\'e duality groups} is~\cite{Davis(2000Poin)}.

The property of  being hyperbolic has a lot of consequences:

\begin{theorem}[Properties of hyperbolic groups] 
\label{the:hyperbolic_groups} \
\begin{enumerate}

\item \label{the:hyperbolic_groups:quasiisometry} \emph{Geometric:}\\
The property ``hyperbolic'' is geometric;

\item \label{the:hyperbolic_groups:geometric_action}
\emph{Characterization by actions:}\\
A group $G$ is a hyperbolic group if and only if it acts
isometrically, properly and cocompactly on a proper hyperbolic space $X$.
In this case $\partial G$ is homeomorphic to $\partial X$;

\item \label{the:hyperbolic_groups:asymptotic_cones}
\emph{Characterization by asymptotic cones:}\\
A finitely generated group is hyperbolic if and only if all its asymptotic cones 
are $\IR$-trees. A finitely presented group is hyperbolic if and only if
one (and hence all) asymptotic cones are $\IR$-trees;

\item \emph{Presentations:}
\label{the:hyperbolic_groups:presentations}\
\begin{enumerate}

\item \label{the:hyperbolic_groups:Dehn_presentation}
A finitely generated group is hyperbolic if and only if it possesses a
Dehn presentation;

\item \label{the:hyperbolic_groups:small_cancellation}
Suppose that the finitely presented group $G$ is a small cancellation group 
in the sense that it admits a presentation which satisfies the condition $C'(1/6)$
or which satisfies both the conditions $C'(1/4)$ and $T(4)$ (see~\cite[Definition~3
in Chapter~8 on page 228]{Ghys-Harpe(1990)}). Then $G$ is hyperbolic;

\end{enumerate}

\item \emph{Classifying spaces and finiteness properties:}
\label{the:hyperbolic_groups:Classifying_spaces_and_finiteness_properties}\
\begin{enumerate}

\item \label{the:hyperbolic_groups:eunderbarG}
If $G$ is hyperbolic, then there exists a finite model for 
the universal space for proper
$G$-actions $\eub{G}$;

\item \label{the:hyperbolic_groups:BG_finite_type}
If $G$ is hyperbolic, then there is a model for $BG$ of finite type,
$H_n(G;\IZ)$ is finitely generated as $\IZ$-module
for $n \ge 0$ and $H_n(G;\IQ)$ is trivial for almost all $n \ge 0$;

\item \label{the:hyperbolic_groups:finitely_presented}
If $G$ is hyperbolic, then $G$ is finitely presented;

\item \label{the:hyperbolic_groups:conjugacy_classes_of_finite_subgroups}
Suppose that $G$ is hyperbolic. Then there are only finitely many conjugacy classes
of finite subgroups;

\item \label{the:hyperbolic_groups:finite_BG}
If $G$ is hyperbolic and torsionfree, then there is a finite model for~$BG$,
the abelian group $H_n(G;\IZ)$ is finitely generated 
for $n \ge 0$ and $H_n(G;\IZ)$ is trivial for almost all $n \ge 0$;

\end{enumerate}

\item  \label{the:hyperbolic_groups:subgroups} \emph{Subgroups:}\
\begin{enumerate}

\item \label{the:hyperbolic_groups:normalizers}
Let $C \subseteq G$ be an infinite cyclic subgroup of $G$. Suppose that
$G$ is hyperbolic.  Then $C$ has finite index
in both its centralizer $C_GC$ and its normalizer~$N_GC$.
In particular, $G$ does not contain a subgroup isomorphic to~$\IZ^n$ for $n \ge 2$;

\item \label{the:hyperbolic_groups:amenable_subgroup}
Any subgroup of a hyperbolic group is either virtually cyclic
or contains a free group of rank two as subgroup. In particular,
an amenable subgroup of a hyperbolic group is virtually cyclic;

\item \label{the:hyperbolic_groups:free_subgroups}
Given $r$ elements $g_1$, $g_2$, \ldots , $g_r$ in a hyperbolic group,
then there exists an integer $n \ge 1$, such that $\{g_1^n, g_2^n,
\ldots, g_r^n\}$ generates a free subgroup of rank at most $r$;

\end{enumerate}

\item  \label{the:hyperbolic_groups:torsion_group} \emph{Torsion groups:}\\
Let $G$ be a torsion group, i.e., each element in $G$ has finite order.
Then $G$ is hyperbolic if and only if $G$ is finite;

\item  \label{the:hyperbolic_groups:inheritance} \emph{Inheritance properties:}\

\begin{enumerate}

\item \label{the:hyperbolic_groups:direct_product}
The product $G_1 \times G_2$ of two hyperbolic groups is again hyperbolic
if and only if one of the two groups $G_1$ and $G_2$ is finite;

\item \label{the:hyperbolic_groups:free_product}
The free product of two hyperbolic groups is again hyperbolic;

\end{enumerate}

\item \label{the:hyperbolic_groups:decision_problems} \emph{Decision problems:}\

\begin{enumerate}

\item \label{the:hyperbolic_groups:conjugation_problem}
The word-problem and the conjugation-problem is solvable for a
hyperbolic group;

\item \label{the:hyperbolic_groups:isomorphism_problem}
The isomorphism-problem is solvable for torsionfree hyperbolic groups;

\end{enumerate}

\item \emph{The boundary:} 
\label{the:hyperbolic_groups:boundary_of_a_hyperbolic_group} \

\begin{enumerate}

\item \label{the:hyperbolic_groups:cd(G)_and_dim(partialG)}
Let $G$ be a hyperbolic group that is virtually torsionfree. Then
$$\vcd(G) - 1 = \dim(\partial G),$$
where $\vcd(G)$ is the virtual cohomological dimension of $G$ and 
$\dim(\partial G)$ is the topological dimension of $\partial G$;

\item \label{the:hyperbolic_groups:boundary}
Let $G$ be hyperbolic and infinite and let $n \ge 2$ be an integer. 
Suppose that $\partial G$ contains an 
open subset which is homeomorphic
to $\IR^n$. Then $\partial G$ is homeomorphic to $S^n$;

\item \label{the:hyperbolic_groups:Fuchsian_groups}
Let $G$ be hyperbolic. Then $\partial G$ is homeomorphic to $S^1$ if and only if
$G$ is a Fuchsian group;

\item \label{the:hyperbolic_groups:Poincare_duality}
A torsionfree hyperbolic group $G$ is a Poincar\'e duality group of dimension~$n$
if and only if $\partial G$ has the integral \v{C}ech cohomology of
$S^{n-1}$;

\item \label{the:hyperbolic_groups:Poincare_duality_in_dim_2}
A torsionfree hyperbolic group $G$ is a Poincar\'e duality group of
dimension~$3$ if and only if $\partial G$ is homeomorphic to $S^2$;

\end{enumerate}

\item \label{the:hyperbolic_groups:rational_function}
\emph{Rationality:}\\[1mm]
Let $G$ be a hyperbolic group. Let $S$ be a finite set of generators.
For the integer $n \ge 0$ let $\sigma(n)$ be the number of elements 
$g \in G$ with $d_S(g,e) = n$;

Then the formal power series $\sum_{n=0}^{\infty} \sigma(n) \cdot t^n$
is a rational function.

The same is true if one replaces $\sigma(n)$ by
the number $\beta(n)$ of elements  $g \in G$ with $d_S(g,e) \le  n$;  

\item \label{the:hyperbolic_groups:further_group_theoretic_properties}
\emph{Further group theoretic properties:}\

\begin{enumerate}

\item \label{the:hyperbolic_groups:weak_amenable}
A hyperbolic group is weakly amenable in the sense 
of Cowling-Haa\-ge\-rup~\cite{Cowling+Haagerup(1989)};

\item \label{the:hyperbolic_groups:finite-asymptotic_dimension}
A hyperbolic group has finite asymptotic dimension;

\item \label{the:hyperbolic_groups:Hopf_Property}
A finitely generated subgroup $H$ of a torsionfree hyperbolic group is
Hopfian, i.e., every epimorphism $H \to H$ is an isomorphism;

\end{enumerate}
\item \emph{Being hyperbolic is generic:}\\[1mm]
\label{the:hyperbolic_groups:statistic}
In a precise statistical sense almost all finitely presented groups are hyperbolic.

\end{enumerate}
\end{theorem}
\begin{proof}%
\ref{the:hyperbolic_groups:quasiisometry} Quasiisometric groups have 
quasiisometric Cayley graphs and it is not difficult to see that
the property being hyperbolic is a quasiisometry invariant of geodesic spaces.
\\[1mm]%
\ref{the:hyperbolic_groups:geometric_action}
See~\cite[Theorem~2.24]{Kapovich+Benakli(2002)}, \cite{Gromov(1987)}. 
\\[1mm]%
\ref{the:hyperbolic_groups:asymptotic_cones}
See~\cite[Section~1.1]{Olshanskii-Osin-Sapir(2007)}.
\\[1mm]%
\ref{the:hyperbolic_groups:Dehn_presentation}
See~\cite[Theorem~2.6 in Chapter~III.$\Gamma$ on page~450]{Bridson-Haefliger(1999)}.
\\[1mm]%
\ref{the:hyperbolic_groups:small_cancellation}
See~\cite[Theorem~36 in Chapter~8 on page 254]{Ghys-Harpe(1990)}.
\\[1mm]%
\ref{the:hyperbolic_groups:eunderbarG}
One can assign to a hyperbolic  group its Rips complex for a certain parameter.
If this parameter is chosen large enough, then the Rips complex is a model
for $\eub{G}$ (see~\cite{Meintrup-Schick(2002)}). The Rips complex
is known to be a $G$-$CW$-complex which is finite or, equivalently, cocompact.
\\[1mm]%
\ref{the:hyperbolic_groups:BG_finite_type} 
This follows from
assertion~\ref{the:hyperbolic_groups:eunderbarG} (see~\cite[Theorem~4.2]{Lueck(2000a)}).
\\[1mm]%
\ref{the:hyperbolic_groups:finitely_presented}
This follows from assertion~\ref{the:hyperbolic_groups:BG_finite_type}.
\\[1mm]%
\ref{the:hyperbolic_groups:conjugacy_classes_of_finite_subgroups}
This follows from
assertion~\ref{the:hyperbolic_groups:eunderbarG} (see~\cite[Theorem~4.2]{Lueck(2000a)}).
\\[1mm]%
\ref{the:hyperbolic_groups:finite_BG}
This follows from
assertion~\ref{the:hyperbolic_groups:eunderbarG}.
\\[1mm]%
\ref{the:hyperbolic_groups:normalizers} 
See~\cite[Corollary~3.10 in Chapter~III.$\Gamma$ on page 462]{Bridson-Haefliger(1999)}.
\\[1mm]%
\ref{the:hyperbolic_groups:amenable_subgroup}
This follows from~\cite[Theorem~37 in Chapter~8 on page~154]{Ghys-Harpe(1990)}
and the fact that an amenable group cannot contain a free group of rank 2 as subgroup.
\\[1mm]%
\ref{the:hyperbolic_groups:free_subgroups}
See~\cite[Proposition~3.20 in Chapter~III.$\Gamma$ on page 467]{Bridson-Haefliger(1999)}.
\\[1mm]%
\ref{the:hyperbolic_groups:torsion_group}
See~\cite[Proposition~2.22 in Chapter~III.$\Gamma$ on page~458]{Bridson-Haefliger(1999)}.
\\[1mm]%
\ref{the:hyperbolic_groups:direct_product}
This follows from assertions~\ref{the:hyperbolic_groups:amenable_subgroup} 
and~\ref{the:hyperbolic_groups:torsion_group}.
\\[1mm]%
\ref{the:hyperbolic_groups:free_product}
See~\cite[Exercise~34 in Chapter~1 on page~19]{Ghys-Harpe(1990)}.
\\[1mm]%
\ref{the:hyperbolic_groups:conjugation_problem}
See~\cite[Theorem~2.8 in Chapter~III.$\Gamma$ on page~451]{Bridson-Haefliger(1999)}.
\\[1mm]%
\ref{the:hyperbolic_groups:isomorphism_problem}
See~\cite{Sela(1995)}.
\\[1mm]%
\ref{the:hyperbolic_groups:cd(G)_and_dim(partialG)}
See~\cite[Corollary~1.4~(e)]{Bestvina-Mess(1991)}.
\\[1mm]%
\ref{the:hyperbolic_groups:boundary}
See~\cite[Theorem~4.4]{Kapovich+Benakli(2002)}.
\\[1mm]%
\ref{the:hyperbolic_groups:Fuchsian_groups}
See~\cite{Casson-Jungreis(1994)}, \cite{Freden(1995)},
\cite{Gabai(1991)}.
\\[1mm]%
\ref{the:hyperbolic_groups:Poincare_duality}
See~\cite[Corollary~1.3]{Bestvina-Mess(1991)}.
\\[1mm]%
\ref{the:hyperbolic_groups:Poincare_duality_in_dim_2}
This follows from
assertion~\ref{the:hyperbolic_groups:Poincare_duality}.
See~\cite[Corollary~6.3]{Davis(2000Poin)}.
\\[1mm]%
\ref{the:hyperbolic_groups:rational_function}
See~\cite[Theorem~2.21 in Chapter~III.$\Gamma$ on page 457]{Bridson-Haefliger(1999)}.
\\[1mm]%
\ref{the:hyperbolic_groups:weak_amenable}
See~\cite{Ozawa(2007)}.
\\[1mm]%
\ref{the:hyperbolic_groups:finite-asymptotic_dimension} 
See~\cite{Roe(2005asym)}.
\\[1mm]%
\ref{the:hyperbolic_groups:Hopf_Property} 
See~\cite{Bumagina(2004)}.
\\[1mm]%
\ref{the:hyperbolic_groups:statistic}
See~\cite{Olshanskii(1992)}.
\end{proof}

\begin{remark}[The boundary of a hyperbolic group]
\label{rem:boundary_of_a_hyperbolic_group}
The boundary $\partial X$ of a hyperbolic space and in particular the boundary
$\partial G$ of a hyperbolic group $G$ are metrizable. 
Any compact metric space can be realized as the boundary of a hyperbolic space.
However, not every compact metrizable space can occur 
as the boundary of a hyperbolic group.
Namely, exactly one of the following three cases occurs:
\begin{enumerate}

\item $G$ is finite and $\partial G$ is empty;

\item $G$ is infinite virtually cyclic and $\partial G$ consists of two points;

\item $G$ contains a free group of rank two as subgroup and $\partial
      G$ is an 
infinite perfect, (i.e., without isolated points) compact metric space.

\end{enumerate}

The metric structure on $\partial X$ for a hyperbolic space $X$ is not canonical. 
One can actually equip $\partial X$ with the structure of
a \emph{visual metric} (see~\cite[Definition~3.20 on page~343]{Bridson-Haefliger(1999)}).
Again the structure of a space with a visual metric
is not canonical, not even for $\partial G$ of a hyperbolic group $G$.
However, the induced \emph{quasiconformal structure} and the induced
\emph{quasi-M\"obius structure} associated to some visual metric on $\partial G$ 
of a hyperbolic group $G$ are canonical, i.e., independent of the choice
of a visual metric. 

These structures are quasiisometry invariants. Namely,
a quasiisometry of finitely generated hyperbolic groups $G_1 \to G_2$ 
(with respect to some choice of finite sets of generators) induces
a homeomorphism $\partial G_1 \to \partial G_2$ 
which is 
quasiconformal and quasi-M\"obius  homeomorphism with respect to any
visual metric.
The converse is also true in the sense that a
homeomorphism $\partial G_1 \to \partial G_2$, which is a quasi-M\"obius 
equivalence or a quasiconformal homeomorphism, comes from a quasiisometry $G_1 \to G_2$.
(see~\cite{Boudon(1995)}, \cite[Section~3]{Kapovich+Benakli(2002)}, \cite{Paulin(1996)}).

The induced action of $G$ on the boundary $\partial G$ is also an important 
invariant of $G$.
\end{remark}

For more information about the boundary of a hyperbolic group we refer
for instance to~\cite{Kapovich+Benakli(2002)}. 

We mention the following result whose proof will appear in a forthcoming paper
by Bartels, L\"uck and  Weinberger~\cite{Bartels-Lueck-Weinberger(2008)}.

\begin{theorem}[High-dimensional spheres as boundary]
\label{the:hyperbolic_groups_and_manifolds}
Let $G$ be a torsionfree hyperbolic group and let $n$ be an integer  $\ge 6$. Then:
\begin{enumerate}

\item \label{the:hyperbolic_groups_and_manifolds:sphere}
The following statements are equivalent:

\begin{enumerate}

\item The boundary $\partial G$ is homeomorphic to $S^{n-1}$;

\item There is a closed aspherical topological manifold $M$ such that $G
\cong \pi_1(M)$, its universal covering $\widetilde{M}$ is homeomorphic to $\IR^n$
and the compactification of $\widetilde{M}$ by $\partial G$ is homeomorphic
to $D^n$.

\end{enumerate}

\item \label{the:hyperbolic_groups_and_manifolds:Cech}
The following statements are equivalent:

\begin{enumerate}

\item The boundary $\partial G$  has the integral \v{C}ech cohomology of
$S^{n-1}$;

\item There is a closed aspherical ANR-homology manifold $M$ with $G
\cong \pi_1(M)$.

\end{enumerate}

\item \label{the:hyperbolic_groups_and_manifolds:uniqueness} 
Let $M$ and $N$ be two aspherical closed $n$-dimensional manifolds together with
isomorphisms $\phi_M \colon \pi_1(M) \xrightarrow{\cong} G$ and 
$\phi_N \colon \pi_1(N) \xrightarrow{\cong} G$. Then there exists
a homeomorphism $f \colon M \to N$ such that $\pi_1(f)$ agrees
with $\phi_N^{-1} \circ \phi_M$ (up to inner automorphisms).

\end{enumerate}
\end{theorem}

\begin{remark}[Algorithm for the homeomorphism problem]
\label{rem:algorithm_for_the_homeomorphism_problem}
By unpublished work of Bartels and L\"uck~\cite{Bartels-Lueck(2008Borel)} 
on the Borel Conjecture for hyperbolic groups
two closed aspherical manifolds with hyperbolic fundamental groups
and dimension $n \ge 5$ are homeomorphic if and only if their fundamental 
groups are isomorphic. Combining this with the result of Sela~\cite{Sela(1995)}
stated in Theorem~\ref{the:hyperbolic_groups}~%
\ref{the:hyperbolic_groups:isomorphism_problem} shows
for any integer $n \ge 5$: There exists an algorithm which takes as
input two closed aspherical $n$-dimensional 
manifolds with hyperbolic fundamental groups
and which (after a finite amount of time) will stop and 
answers \emph{yes} or \emph{no} according to whether
or not the manifolds are homeomorphic.

The following is already pointed out in~\cite[page~459]{Bridson-Haefliger(1999)}:
There is a technical problem here with how the closed aspherical manifolds are given. 
They must be given by a finite amount of information 
(from which one can read off a presentation
of the fundamental group). 
\end{remark}

\begin{remark}[Lacunary groups]\label{rem:lacunary_groups}
Olshanskii-Osin-Sapir~\cite{Olshanskii+Sapir(2006)}  introduced the notion of
a \emph{lacunary} group as a finitely generated 
group one of whose asymptotic cones is an $\IR$-tree.
They show that such a group can always be obtained as a colimit of
a directed system  of hyperbolic groups
$G_1 \to G_2 \to G_3 \to \cdots,$
where the structure maps are epimorphisms of hyperbolic groups
with certain additional properties.
A finitely presented lacunary group is hyperbolic. The class of lacunary
groups is very large and contains some examples with unusual properties, e.g.,
certain infinite torsionfree groups whose proper subgroups are all cyclic and
infinite torsion-groups whose proper subgroups are all of order $p$ for some fixed prime 
number $p$.

\end{remark}

\begin{remark}
\label{rem:hyperbolic_groups_yield_exotic_groups}
Colimits of directed systems of
hyperbolic groups  which come from adding more and more relations
have been used to construct exotic groups. 
Other constructions come from random groups (see~\cite{Gromov(2003)}).
Here are some examples:
\begin{enumerate}

\item Let $G$ be a torsionfree hyperbolic group which is not virtually cyclic.
      Then there exists a quotient of $G$ which is an infinite torsiongroup
      whose proper subgroups are all finite (or cyclic) (See~\cite{Olshanskii(1993)});

\item There are hyperbolic groups which do have  Kazhdan's property
      (T) (see Zuk~\cite{Zuk(2003)});

\item There exist groups with expanders. They play a role in the construction
       of counterexamples to the Baum-Connes Conjecture with coefficients
       due to Higson, Lafforgue and Skandalis~\cite{Higson-Lafforgue-Skandalis(2002)}.

\end{enumerate}
\end{remark}

\begin{remark}[Exotic aspherical manifolds]
For every $n \ge 5$ there exists an example of a closed aspherical
topological manifold $M$ of dimension $n$  
that is a piecewise flat, non-positively curved polyhedron such that the
universal covering $\widetilde{M}$ is not homeomorphic to $\IR^n$
(see~\cite[Theorem~5b.1 on page~383]{Davis-Januszkiewicz(1991)}).
This manifold is not homeomorphic to a closed smooth manifold with Riemannian metric of non-positive
sectional curvature by Hadamard's Theorem.
There is a variation of this construction that uses the 
strict hyperbolization of Charney-Davis~\cite{Charney-Davis(1995)} 
and produces closed aspherical manifolds whose universal cover is 
not homeomorphic to 
Euclidean space and whose fundamental group is hyperbolic.

There exists a strictly negatively curved polyhedron $N$ of dimension $5$
whose fundamental group is hyperbolic, which is homeomorphic to a
closed aspherical smooth manifold and 
whose universal covering is homeomorphic to $\IR^n$, but 
the ideal boundary of its universal covering, which is homeomorphic to
$\partial G$, is not homeomorphic to $S^{n-1}$
(see~\cite[Theorem~5c.1 on page~384]{Davis-Januszkiewicz(1991)}).
Notice $N$ is not homeomorphic to a closed smooth Riemannian
manifold with negative sectional curvature.
\end{remark}

\begin{remark}[Cohomological characterization of hyperbolic groups]
\label{rem:Cohomological–characterizations_of_hyperbolic_groups}
There exist also characterizations of the property hyperbolic in terms
of cohomology. A finitely presented group $G$ is hyperbolic if and
only if $H^{(1)}_1(G,\IR) = \overline{H}^{(1)}_1(G,\IR) = 0$ holds for
the first $l^1$-homology and the first reduced $l^1$-homology
(see~\cite{Allcock+Gersten(1999)}).
For a characterization in terms of bounded cohomology we refer to~\cite{Mineyev(2002boundcoh)}. 
\end{remark}


\typeout{--------------------   Section 8 ---------------------------------------}

\section{CAT(0)-groups}
\label{sec:CAT(0)-groups}

\begin{definition}[$\CAT(0)$-group]\label{def:CAT(0)-group}
A group is called \emph{$\CAT(0)$-group} if it admits an isometric proper
cocompact action on some $\CAT(0)$-space.
\end{definition}

\begin{theorem}[Properties of $\CAT(0)$-groups]
\label{the:CAT(0)-groups}
\
\begin{enumerate}

\item \label{the:CAT(0)-groups:Classifying_spaces_and_finiteness_properties:}
\emph{Classifying spaces and finiteness properties:}\

\begin{enumerate}

\item \label{the:CAT(0)-groups:eunderbarG}
If $G$ is a $\CAT(0)$-group, then there exists a finite model for 
the universal space of proper
$G$-actions $\eub{G}$ (see Definition~\ref{def:classifying_space_for_proper_actions});

\item \label{the:CAT(0)-groups:BG_of_finite_type}
If $G$ is a $\CAT(0)$-group, then there is a model for $BG$ of finite type,
$H_n(G;\IZ)$ is finitely generated as $\IZ$-module
for $n \ge 0$ and $H_n(G;\IQ)$ is trivial for almost all $n \ge 0$;

\item \label{the:CAT(0)-groups:finitely_presented}
If $G$ is a $\CAT(0)$-group, then $G$ is finitely presented;

\item \label{the:CAT(0)-groups:conjugacy_classes_of_finite_subgroups}
Suppose that $G$ is a $\CAT(0)$-group. Then there are only finitely many 
conjugacy classes of finite subgroups of $G$;

\item \label{the:CAT(0)-groups:finite_BG}
If $G$ is a torsionfree $\CAT(0)$-group, then there is a finite model
for $BG$, the abelian group $H_n(G;\IZ)$ is finitely generated
for $n \ge 0$ and $H_n(G;\IZ)$ is trivial for almost all $n \ge 0$;

\end{enumerate}

\item\label{the:CAT(0)-groups:solvable}  \emph{Solvable subgroups:} \\
Every solvable subgroup of a $\CAT(0)$-group is virtually $\IZ^n$;

\item\label{the:CAT(0)-groups:Inheritance_properties} \emph{Inheritance properties:}\
\begin{enumerate}

\item \label{the:CAT(0)-groups:direct_products}
The direct product of two $\CAT(0)$-groups is again a $\CAT(0)$-group;

\item \label{the:CAT(0)-groups:amalgamation}
The free product with amalgamation along a virtually cyclic subgroup
of two $\CAT(0)$-groups is again a $\CAT(0)$-group;

\item \label{the:CAT(0)-groups:HNN}
The $\HNN$-extension of a $\CAT(0)$-group along a finite group is again
a $\CAT(0)$-group;

\end{enumerate}

\item \label{the:CAT(0)-groups:exmaples} \emph{Examples:}\
\begin{enumerate}

\item \label{the:CAT(0)-groups:limit_groups}
Limit groups in the sense of Sela are $\CAT(0)$-groups;

\item \label{the:CAT(0)-groups:Coxeter_groups}
Coxeter groups are $\CAT(0)$-groups;

\item \label{the:CAT(0)-groups:special_Artin_groups}
Three-dimensional FC Artin groups are $\CAT(0)$-groups;

\end{enumerate}

\item \label{the:CAT(0)-groups:conjugation_problem} \emph{Decision problems:}\\
The word-problem and the conjugation-problem are solvable for a $\CAT(0)$-group;

\item \label{the:CAT(0)-groups:versus_hyperbolic} \emph{Hyperbolic:}\\
Let $G$ act isometrically, properly and cocompactly on the $\CAT(0)$-space $X$.
Then $G$ is hyperbolic if and only if $X$ does not contain an isometrically
embedded copy of a Euclidean plane;

\item \label{the:CAT(0)-groups:weak_hyperbolization}
\emph{Weak Hyperbolization Theorem:}\\
Let $G$ be a three-dimensional Poincar\'e duality group. Suppose
that in addition that $G$ is a $\CAT(0)$-group. Then $G$ satisfies the Weak
Hyperbolization Conjecture, i.e., either $G$ contains $\IZ^2$ or $G$ is hyperbolic.

\end{enumerate}
\end{theorem}
\begin{proof}%
\ref{the:CAT(0)-groups:eunderbarG}
Let $X$ be a $\CAT(0)$-space on which $G$ acts properly, isometrically and cocompactly.
Then it is easy to show that $X$ is a model for $\underline{J}G$ for the numerable
version of the classifying space for proper $G$-actions. (Notice that $X$ is not necessarily
a $CW$-complex. But this implies that there is a cocompact model for $\eub{G}$.
Details will appear in~\cite{Lueck(2009catevcyc)}.
\\[1mm]%
\ref{the:CAT(0)-groups:BG_of_finite_type}
This follows from
assertion~\ref{the:CAT(0)-groups:eunderbarG} (see~\cite[Theorem~4.2]{Lueck(2000a)}).
\\[1mm]%
\ref{the:CAT(0)-groups:finitely_presented}
This follows from assertion~\ref{the:CAT(0)-groups:BG_of_finite_type}.
\\[1mm]%
\ref{the:CAT(0)-groups:conjugacy_classes_of_finite_subgroups}
This follows from
assertion~\ref{the:CAT(0)-groups:eunderbarG} (see~\cite[Theorem~4.2]{Lueck(2000a)}).
\\[1mm]%
\ref{the:CAT(0)-groups:finite_BG}
This follows from assertion~\ref{the:CAT(0)-groups:eunderbarG}.\\[1mm]%
\ref{the:CAT(0)-groups:solvable}
See~\cite[Theorem~1.1 in Chapter III.$\Gamma$ on page~439]{Bridson-Haefliger(1999)}.
\\[1mm]%
\ref{the:CAT(0)-groups:direct_products}
See~\cite[Theorem~1.1 in Chapter III.$\Gamma$ on page~439]{Bridson-Haefliger(1999)}.
\\[1mm]%
\ref{the:CAT(0)-groups:amalgamation}
See~\cite[Theorem~1.1 in Chapter III.$\Gamma$ on page~439]{Bridson-Haefliger(1999)}.
\\[1mm]%
\ref{the:CAT(0)-groups:HNN}
See~\cite[Theorem~1.1 in Chapter III.$\Gamma$ on page~439]{Bridson-Haefliger(1999)}.
\\[1mm]%
\ref{the:CAT(0)-groups:limit_groups}
\cite{Alibegovic+Bestvina(2006)}.
\\[1mm]%
\ref{the:CAT(0)-groups:Coxeter_groups}
This is a result due to Moussong. 
See~\cite[Theorem~12.3.3 on page~235]{Davis(2008cox)},
\cite{Moussong(1987)}.
\\[1mm]%
\ref{the:CAT(0)-groups:special_Artin_groups}
See~\cite{Bell(2005)}.
\\[1mm]%
\ref{the:CAT(0)-groups:versus_hyperbolic}
See~\cite[Theorem~3.1 in Chapter III.$\Gamma$ on page~459]{Bridson-Haefliger(1999)}.
\\[1mm]%
\ref{the:CAT(0)-groups:weak_hyperbolization} 
See~\cite[Theorem~2]{Kapovich+Kleiner(2007)}.
\end{proof}

Interesting results about $\CAT(0)$-groups and $\CAT(0)$-lattices including rigidity statements
have been proved by Caprace and Monod~\cite{Caprace-Monod(2008)}.


\typeout{--------------------   Section 9 ---------------------------------------}

\section{Classifying spaces for proper actions}
\label{sec:Classifying_spaces_for_proper_actions}

Very often information or basic properties of groups are reflected
in interesting actions of the group. In this context the notion of a classifying space 
for proper $G$-actions is important. This notion and the more general 
notion of a classifying space for a family of subgroups was introduced by 
tom Dieck (see~\cite{Dieck(1972)}, \cite[I.6]{Dieck(1987)}).

A \emph{$G$-$CW$-complex} $X$ is a $CW$-complex with a 
$G$-action such that for every open cell
$e$ and every $g \in G$ with $g \cdot e = e$ we have $gx=x$ 
for every $g \in G$ and $x \in e$. The barycentric subdivision
of a simplicial complex with simplicial $G$-action is a $G$-$CW$-complex.
A $G$-$CW$-complex $X$ is \emph{proper}
if and only if all its isotropy groups are finite
(see~\cite[Theorem 1.23]{Lueck(1989)}).

\begin{definition}[Classifying space for proper actions]
\label{def:classifying_space_for_proper_actions}
Let $G$ be a group. A model for the \emph{classifying space of proper
$G$-actions} is a proper $G$-$CW$-complex $\eub{G}$ such that 
$\eub{G}^H$ is contractible for all finite subgroups $H \subseteq G$.
\end{definition}

\begin{theorem} [Homotopy characterization of $\eub{G}$]
\label{the:homotopy_characterization_of_eub(G)} \ 
\begin{enumerate}

\item \label{the:homotopy_characterization_of_eub(G):existence}
There exists a model for $\eub{G}$;

\item \label{the:homotopy_characterization_of_eub(G):characterization}
A $G$-$CW$-complex $Y$ is a model for $\eub{G}$ if and only if 
for every proper $G$-$CW$-complex $X$ there is up to $G$-homotopy
precisely one $G$-map $X \to Y$. In particular any two models
for the classifying space for proper $G$-actions are $G$-homotopy equivalent.
 
\end{enumerate}

\end{theorem}
\begin{proof} See for instance~\cite[Theorem~1.9 on page~275]{Lueck(2005s)}.
\end{proof}

If $G$ is torsionfree, then a model for $\eub{G}$ is a model for $EG$, i.e., 
the total space of the universal $G$-principal bundle $G \to EG \to BG$.
A group $G$ is finite if and only if $G/G$ is a model for $\eub{G}$.

Some prominent groups come with prominent actions on prominent spaces.
Often it turns out that these are models for the classifying space for proper $G$-actions.
Here we give a  list of examples. More explanations and references can be found
in the survey article~\cite{Lueck(2005s)}.

\begin{itemize}

\item Discrete subgroups of almost connected Lie groups\\
Let $L$ be a Lie group with finitely many path components. 
Let $K \subseteq L$ be any maximal compact subgroup,
which is unique up to conjugation. Let $G \subseteq L$ be a discrete subgroup.
Then $L/K$ is diffeomorphic to $\IR^n$ and becomes with the obvious left $G$-action a model
for $\eub{G}$.

\item Hyperbolic groups and the Rips complex\\
Let $G$ be a hyperbolic group. Let $P_d(G)$ be the Rips complex.
Then $P_d(G)$ is a model for $\eub{G}$ if $d$ is chosen large enough.

\item Proper isometric actions on simply connected complete Riemannian
 manifolds with non-positive sectional curvature\\
Suppose that $G$ acts isometrically and properly on a simply connected
complete Riemannian manifold $M$ with non-positive sectional curvature.
Then $M$ is a model for $\eub{G}$;

\item Proper actions on trees\\
Let $T$ be  a tree. Suppose that $G$ acts on $T$ by tree automorphisms
without inversion such that all isotropy groups are finite.
Then $T$ is a model for~$\eub{G}$;

\item Arithmetic groups and the Borel-Serre compactification\\
Let $G(\IR)$ be the $\IR$-points of a semisimple
$\IQ$-group $G(\IQ)$  and let $K\subseteq G(\IR)$ be a maximal compact subgroup.
If $A \subseteq G(\IQ)$ is an arithmetic group, then $G(\IR)/K$ with the left 
$A$-action is a model for $\eub{A}$. The $A$-space $G(\IR)/K$ is not necessarily cocompact.
However, the Borel-Serre completion  of  $G(\IR)/K$ is a finite $A$-$CW$-model 
for $\eub{A}$; 

\item Mapping class groups and Teichm\"uller space\\
Let $\Gamma^s_{g,r}$be the \emph{mapping class group} of an orientable compact surface
$F^s_{g,r}$ of genus $g$ with $s$ punctures and $r$ boundary components.
This is the group of isotopy classes of orientation preserving
self-diffeomorphisms $F^s_{g,r} \to F^s_{g,r}$ that  preserve the punctures 
individually and restrict to the
identity on the boundary. We require that the isotopies
leave the boundary pointwise fixed. We will always assume
that $2g +s +r > 2$, or, equivalently, that the Euler characteristic
of the punctured surface $F^s_{g,r}$ is negative. Then the associated
\emph{Teich\-m\"ul\-ler space} $\calt^s_{g,r}$
is a model for $\eub{\Gamma^s_{g,r}}$;

\item $\Out(F_n)$ and outer space\\
Let $F_n$ be the free group of rank $n$. Denote by 
$\Out(F_n)$ the group of outer automorphisms
of $F_n$. Culler and Vogtmann~\cite{Culler-Vogtmann(1986)}, 
\cite{Vogtmann(2002)}
have constructed a space $X_n$ called \emph{outer space},
on which $\Out(F_n)$ acts with finite isotropy groups. 
It is a model for $\eub{\Out(F_n)}$.

The space $X_n$ contains a \emph{spine}
$K_n$ which is an $\Out(F_n)$-equivariant deformation retract.
This space $K_n$ is a simplicial complex of dimension $(2n-3)$
on which the $\Out(F_n)$-action is by simplicial automorphisms
and cocompact.  Hence the barycentric subdivision of $K_n$ is a finite
$(2n-3)$-dimensional model of $\underline{E}\Out(F_n)$;

\item One-relator groups\\
Let $G$ be a one-relator group. Let $G = \langle (q_i)_{i \in I} \mid r \rangle$ 
be a presentation with one relation. 
There is up to conjugacy one maximal finite subgroup~$C$ which turns out to be cyclic.
Let $p \colon \ast_{i \in I} \IZ \to G$ be the epimorphism from the free group
generated by the set $I$ to $G$ that sends the generator $i \in I$
to~$q_i$.  Let $Y \to \bigvee_{i \in I} S^1$ be the $G$-covering associated to the epimorphism
$p$. There is a $1$-dimensional  unitary $C$-representation $V$
and a $C$-map $f \colon SV \to \res_G^C Y$ such that the following is true:
The induced action on the unit sphere $SV$ is free. 
If we equip $SV$ and $DV$ with the obvious
$C$-$CW$-complex structures, the $C$-map $f$ can be chosen to be cellular 
and we obtain a $G$-$CW$-model for $\eub{G}$
by the $G$-pushout
$$\xycomsquareminus{G \times_C SV}{\overline{f}}{Y}{}{}{G \times_C DV}{}{\eub{G}}$$
where $\overline{f}$ sends $(g,x)$ to $gf(x)$. 
Thus we get a $2$-dimensional $G$-$CW$-model for $\underline{E}G$ such that
$\eub{G}$ is obtained from $G/C$ for a maximal finite cyclic
subgroup $C \subseteq G$ by attaching free cells of dimensions $\le 2$
and the $CW$-complex structure on the quotient $G\backslash \underline{E}G$ 
has precisely one $0$-cell, precisely one $2$-cell and as many $1$-cells as there are
elements in $I$.

\end{itemize}

\begin{remark}[Isomorphism Conjectures]
\label{rem:Isomorphism_Conjectures}
The space $\eub{G}$ and its version for the family of virtually cyclic
subgroups play an important role in the formulation of the Isomorphism
Conjectures for $K$- and $L$-theory of group rings and reduced group $C^*$-algebras
or Banach algebras due to Farrell-Jones
(see~\cite[1.6 on page 257]{Farrell-Jones(1993a)}), Baum-Connes 
(see~\cite[Conjecture 3.15 on page 254]{Baum-Connes-Higson(1994)}) and Bost.
Methods and results from
geometric group theory enter the proofs of these conjectures for certain classes
of groups. A survey on these conjectures, their status and the methods of proof
can be found for instance in~\cite{Lueck-Reich(2005)}.
\end{remark}

\begin{remark}[Small models]
\label{rem:small_models} 
As one can ask whether there are small models for $BG$ 
(or, equivalently, for the $G$-$CW$-complex $EG$) such as 
finite models, models of finite type or finite-dimensional models, the same question
is interesting for the $G$-$CW$-complex $\eub{G}$ and has been studied for instance 
in~\cite{Kropholler-Mislin(1998)}, \cite{Lueck(2000a)}, \cite{Lueck-Meintrup(2000)}.

Although there are often nice small models for $\eub{G}$, these spaces can be arbitrarily
complicated. Namely, for any $CW$-complex $X$ there exists a group $G$ such that
$G\backslash \eub{G}$ and $X$ are homotopy equivalent
(see~\cite{Leary-Nucinkis(2001a)}). There can also be dramatic changes
in the complexity and size of $\eub{G}$ if one passes from $\eub{H}$ to $\eub{G}$
for a subgroup $H \subseteq G$ of finite index (see~\cite{Leary-Nucinkis(2003)}).
\end{remark}

\begin{remark}[Compactifications of $\eub{G}$]
It is very important to find appropriate compactifications of
$\eub{G}$. Finding the right one which is  ``small at infinity'' leads
to injectivity results concerning the Isomorphism Conjectures  
(see for instance~\cite{Carlsson-Pedersen(1995a)},
\cite{Rosenthal(2004)}, \cite{Rosenthal-Schuetz(2005)}). We have seen for a hyperbolic
group that its boundary yields a powerful compactification of the associated
Rips complex. A $\CAT(0)$-space comes with a natural compactification by adding 
its boundary. There is a whole theory of compactifications of the Teichm\"uller space.
For arithmetic groups the Borel-Serre compactification is crucial.
\end{remark}

\begin{remark}[Computations]
\label{rem:computations}
A good understanding of the spaces $\eub{G}$ can be used to make 
explicit computations of the homology or topological $K$-theory  $H_*(BG)$ and $K_*(BG)$
or various $K$- and $L$-groups such as $K_*(RG)$, $L_*(RG)$ and $K_*(C^*_r(G))$.
See for instance~\cite{Lueck(2002b)}, \cite{Lueck(2002d)},
\cite{Lueck(2005heis)},  \cite{Lueck(2007)},
\cite{Lueck-Stamm(2000)}.
\end{remark}


\typeout{--------------------   Section 10 ---------------------------------------}

\section{Measurable group theory}
\label{sec:Measurable_group_theory}

Gromov~\cite[0.2.$C_2'$ on page~6]{Gromov(1993)}
(see also~\cite[Exercise 35 in~IV.B on page~98]{delaHarpe(2000)}
or~\cite[Theorem~2.1.2]{Shalom(2004)})
observed that the notion of quasiisometry can be reformulated as follows.

\begin{lemma} \label{lem.equivalent_reformulation_of_quasiisometry}
Two finitely generated groups $G_1$ and $G_2$ are quasiisometric if and only if
there exists a locally compact space on which $G_1$ and $G_2$ act properly
and cocompactly and the actions commute.
\end{lemma}

This led Gromov to the following measure theoretic version
(see~\cite[0.5E]{Gromov(1993)}, \cite{Furman(1999a)} and~\cite{Furman(1999b)}).
A \emph{Polish space} is a separable  topological space
which is metrizable by a complete metric. A measurable space  is called a 
\emph{standard Borel space} 
if it is isomorphic to a Polish space with its standard Borel $\sigma$-algebra.
Let $\Omega$ be a standard Borel space with a Borel measure $\mu$. 
Let $G$ act on $\Omega$ by Borel automorphisms. A 
\emph{measure theoretic fundamental domain}
for the $G$-action is a Borel subset $X \subseteq \Omega$ such that
$\mu(g\cdot X \cap X) = 0$ for every $g \in G, g \not= 1$ 
and $\mu(\Omega - G \cdot X) = 0$ hold.

\begin{definition}[Measure equivalence]
\label{def:measure_equivalence}
Two countable groups $G$ and $H$ are called
\emph{measure equivalent} if there is a standard Borel 
space $\Omega$ with a non-zero Borel measure
on which $G$ and $H$ act by measure-preserving Borel automorphisms
such that the actions commute and the
actions of both $G$ and $H$ admit  finite measure fundamental domains.
\end{definition}

The actions appearing in Definition~\ref{def:measure_equivalence} are automatically
\emph{essentially free},
i.e., the stabilizer of almost every point is trivial,
because of the existence of the measure fundamental domains.
Measure equivalence defines an equivalence relation on countable groups
(see~\cite[Section~2]{Furman(1999a)}).

\begin{remark}[Lattices]
Let $\Gamma$ and $\Lambda$ 
be two lattices in the locally compact second countable topological 
group $G$, i.e., discrete subgroups with finite covolume with respect to a Haar measure
on $G$.  Then $\Lambda$  and $\Gamma$ are quasiisometric provided that they are cocompact.
An important feature of measure equivalence is that
$\Lambda$ and $\Gamma$ are measure equivalent without the hypothesis of being
cocompact (see~\cite[0.5.$E_2$]{Gromov(1993)}).
\end{remark}

An action $G \action X$ of a countable group
$G$ is called \emph{standard} if $X$ is a standard Borel space with a probability measure
$\mu$, the group $G$ acts by $\mu$-preserving Borel automorphisms
and the action is essentially free.

\begin{definition}[(Weak) orbit equivalence]
Two standard actions $G \action X$ and $H \action Y$ are called
\emph{weakly orbit equivalent} if there exist
Borel subsets $A \subseteq X$ and $B \subseteq Y$ meeting almost every
orbit  and a Borel isomorphism $f \colon A \to B$ which preserves 
the normalized measures on $A$ and $B$, respectively, and satisfies for almost all $x \in A$
$$f(G\cdot x \cap A) = H \cdot f(x) \cap B.$$
If $A$ and $B$ have full measure in $X$ and $Y$, the two actions are called
\emph{orbit equivalent}.
\end{definition}

The following result is formulated  and proved 
in~\cite[Theorem~3.3]{Furman(1999b)}, where credit is also given to Gromov and Zimmer.

\begin{theorem}[Measure equivalence versus weak orbit equivalence]
\label{the:Measure_equivalence_versus_weak_orbit_equivalence}
Two countable groups $G$ and $H$ are measure equivalent if and only if
there exist standard actions of $G$ and $H$ that are weakly orbit equivalent.
\end{theorem}

The next result is due to Ornstein-Weiss~\cite{Ornstein-Weiss(1980)}.
\begin{theorem} \label{the:Ornstein-Weiss}\
\begin{enumerate}
\item Let $G_1$ and $G_2$ be two infinite countable amenable groups.
Then any two standard actions of $G_1$ and $G_2$ are orbit equivalent;

\item Any infinite amenable group $G$ is measure equivalent to $\IZ$.
\end{enumerate}
\end{theorem}

On the other hand we have the following result due to
Epstein~\cite[Corollary~1.2]{Epstein(2008)}, 
the case of a group with property (T)
has been treated by Hjorth~\cite{Hjorth(2005)} before.

\begin{theorem}\label{the:Epstein}
A countable non-amenable  group admits a continuum of standard
actions which are not pairwise orbit equivalent.
\end{theorem}

The following result is due to
Gaboriau-Popa~\cite{Gaboriau-Popa(2005)}.

\begin{theorem}\label{the:Gaboriau-Popa}
Let $G$ be a non-abelian free group. Then there exists
a continuum of standard
actions $G \action X$  which are pairwise not orbit equivalent
and whose associated von Neumann algebras 
$L^{\infty}(X) \rtimes G$ are pairwise not isomorphic.
\end{theorem}

\begin{remark}[Quasiisometry versus measure equivalence]
\label{rem:Quasiisometry_versus_measurable_equivalence}
In general two finitely presented measure equivalent groups need not  be quasiisometric.
For example $\IZ^n$ and $\IZ^m$ for $n,m \ge 1$ are quasiisometric if and only if
$n = m$ (see Theorem~\ref{the:quasiisometry_invariants})
and they are always measure equivalent 
(see Theorem~\ref{the:Measure_equivalence_versus_weak_orbit_equivalence}).

We mention that property (T) is invariant under measure equivalence 
(see~\cite[Theorem~8.2]{Furman(1999a)}) but is not a quasiisometry
invariant (see Remark~\ref{rem:property_(T)}). 

In general two finitely presented quasiisometric groups need not  be measure equivalent
as the following example shows. If $F_g$ denotes the free group on $g$ generators, then
define $G_n := (F_3 \times F_3) \ast F_n$ for $n \ge 2$.
The groups $G_m$ and $G_n$ are
quasiisometric for $m,n \ge 2$ (see \cite[page~105 in~IV-B.46]{delaHarpe(2000)},
\cite[Theorem 1.5]{Whyte(1999)})
and have finite models for their classifying spaces.
One easily checks that $b_1^{(2)}(G_n) = n$ and $b_2^{(2)}(G_n) = 4$
(see~\cite[Example~1.38 on page~41]{Lueck(2002)}). By the following result of 
Gaboriau~\cite[Theorem 6.3]{Gaboriau(2002a)}  the groups 
$G_n$ and $G_m$ are measure equivalent if and only if $m = n$ holds.
\end{remark}

\begin{theorem}[Measure equivalence and $L^2$-Betti numbers]
\label{the:Gaboriaus_result_on_measure_equivalent_groups}
Let $G_1$ and $G_2$ be two countable groups that are measure equivalent.
Then there is a constant $C > 0$ such that for all $p \ge 0$
$$b_p^{(2)}(G_1)  =  C \cdot b_p^{(2)}(G_2).$$
\end{theorem}

\begin{remark}[Measure equivalence rigidity]
\label{rem:measure_equivalence_rigidity}
In view of Theorem~\ref{the:Ornstein-Weiss} one realizes that measurable
equivalence cannot capture any group theoretic property which can be separated
within the class of amenable groups and is highly non-rigid for amenable groups.
Nevertheless, there is a deep and interesting rigidity theory underlying the 
notion of orbit equivalence.  For information about this topic we refer for instance to the
survey article of Shalom~\cite{Shalom(2005)}. We give as an illustration 
some examples below.
\end{remark}

The next result follows from Furman~\cite[Corollary~B]{Furman(1999b)}
and is stated in the present sharpened form in~\cite[Theorem~3.1]{Shalom(2005)}.

\begin{theorem}\label{the:Furmann}
Fix an odd natural number $n \ge 3$. Consider the obvious standard action of $\mathrm{SL}_n(\IZ)$ on the $n$-torus
$T^n$ equipped with the Lebesgue measure. Suppose that it is orbit equivalent
to a standard action of the group $\Lambda$. Then $\Lambda \cong \mathrm{SL}_n(\IZ)$ and
the orbit equivalence is induced by an isomorphism of actions.
\end{theorem}

The next result is due to Monod-Shalom~\cite[Theorem~1.18]{Monod+Shalom(2006)}
and may be viewed as the measure theoretic definition of a negatively curved group.

\begin{theorem} \label{the:H2b(G,l2(G))_not_0}
The condition that the second bounded cohomology 
$H^2_b(G,l^2(G))$ with coefficients
in $l^2(G)$ does not vanish is an invariant under measure equivalence.
Hyperbolic groups have this property.
\end{theorem}

The next result is taken from~\cite[Corollary~1.11 and Theorem~1.16]{Monod+Shalom(2006)}.
For a countable group $G$ and any probability distribution $\mu$
(different from Dirac) on the interval $[0,1]$, the natural shift action
on $\prod_{G} ([0,1],\mu)$ is called a \emph{Bernoulli $G$-action}.

\begin{theorem}\label{the:Monod-Shalom_Bernoulli_actions}\
\begin{enumerate}

\item Let $G$ be the direct product of two torsionfree groups $G_1$ and $G_2$
with non-trivial  $H^2_b(G_1,l^2(G_1))$ and $H^2_b(G_2,l^2(G_2))$.
If a Bernoulli $G$-action is orbit equivalent to a Bernoulli
$H$-action for some group $H$, then
$G \cong H$ and the actions are isomorphic by a Borel isomorphism
which induces the given orbit equivalence;

\item Let $G_1, \ldots, G_m$ and $H_1, \ldots, H_n$ be torsionfree groups
with non-vanishing $H^2_b(G_i,l^2(G_i))$ and $H^2_b(H_j,l^2(H_j))$.
Suppose that $\prod_{i=1}^m G_i$ and $\prod_{j=1}^n H_j$ are measure equivalent.
Then $m = n$ and for an appropriate permutation $\sigma$ the groups
$G_i$ and $H_{\sigma(i)}$ are measure equivalent for $i = 1,2, \ldots n$.

\end{enumerate}
\end{theorem}

There are many more interesting results in the spirit that the orbit structure of an action 
remembers the group and the action, and relations between orbit equivalence and questions 
about von Neumann algebras and bounded cohomology. In particular, Popa has 
proved spectacular results on fundamental groups of II$_1$-factors and
on rigidity. See for 
instance~\cite{Burger+Monod(2002)},
\cite{Gaboriau(2005)},
\cite{Gaboriau-Popa(2005)},
\cite{Monod+Shalom(2006)},
\cite{Ozawa+Popa(2004)},
\cite{Popa(2004)},
\cite{Popa(2006_betti)},
\cite{Popa(2006_rigidity_bernoulli)},
\cite{Popa(2006_strong_rigidity_actions_I)},
\cite{Popa(2006_strong_rigidity_actions_II},
\cite{Popa(2007)},
\cite {Popa(2007_rigidity)},
\cite{Popa+Vaes(2008)}.


\typeout{--------------------   Section 11 ---------------------------------------}

\section{Some open problems}
\label{sec:Some_open_problems}

Here is a list of interesting open problems.
It reflects some of the interests (and limited knowledge) of the author:


\subsection{Hyperbolic groups}

\begin{enumerate}
\item Is every hyperbolic group virtually torsionfree?

\item Is every hyperbolic group residually finite?

\item Suppose that the space at infinity of a hyperbolic group is
homeomorphic to $S^2$. Does this imply that it acts properly
isometrically and cocompactly on the $3$-dimensional hyperbolic space?

Partial results in this direction have been proved in~\cite{Bonk-Kleiner(2005)}.

\item Has the boundary of a hyperbolic group the integral \v{C}ech cohomology of a sphere
if and only if it occurs as the fundamental group of an  aspherical closed manifold~$M$?

\item Is the boundary $\partial G$ of a hyperbolic group $G$ homeomorphic to $S^n$
if and only if it occurs as the fundamental group of an aspherical closed manifold $M$
whose universal covering $\widetilde{M}$ is homeomorphic to $\IR^n$ and
its compactification $\widetilde{M} \cup \partial G$
by $\partial G$ is homeomorphic to $D^n$?

The answer for the last  problems is yes for $n \ge 6$ 
(see Theorem~\ref{the:hyperbolic_groups_and_manifolds}).

\item Which  topological spaces occur as boundary
      of a hyperbolic group?

\item Is every hyperbolic group  a $\CAT(0)$-group?

\end{enumerate}


\subsection{Isomorphism Conjectures}

\begin{enumerate}

\item Are the Conjectures due to Baum-Connes, Farrell-Jones and Borel true
for the following groups?
\begin{itemize}
\item $\mathrm{SL}_n(\IZ)$ for $n \ge 3$; 
\item Mapping class groups;
\item $\operatorname{Out}(F_n)$.
\end{itemize}

The fact that these conjectures are not known for these groups indicates
that we do not understand enough about the geometry of these groups.
Probably any successful proof will include 
new interesting information about these groups.

\item Are the Conjectures due to Farrell-Jones and Borel true for amenable groups?

The Baum-Connes Conjecture is known for groups with the Haagerup property and hence
in particular for amenable groups
(see~\cite{Higson-Kasparov(2001)}, \cite{Farley(2003)}.)
The Farrell-Jones Conjecture and the Borel Conjecture for these groups
are harder since there one has to take into account all virtually cyclic subgroups
and not only all finite subgroups as in the Baum-Connes setting and one encounters
Nil-phenomena which do not occur in the Baum-Connes setting.

\item Is there a property for groups known such that Isomorphism
Conjectures mentioned above are not known for any group having this property.
If yes, can one use this property to produce counterexamples?

For some time property (T) was thought to be such a property for the Baum-Connes
Conjecture until Lafforgue (see~\cite{Lafforgue(1999)}, \cite{Lafforgue(2001)}) 
proved the Baum-Connes Conjecture 
for certain groups having property (T). 
The counterexamples to the Baum-Connes Conjecture
by Higson-Lafforgue-Skandalis~\cite{Higson-Lafforgue-Skandalis(2002)}
given by groups with expanders have indicated another source of possible counterexamples.
Such groups can be constructed by directed colimits of hyperbolic groups.
However, for directed colimits of hyperbolic groups 
the Farrell-Jones Conjecture and the Borel Conjecture 
in dimension $\ge 5$ are known 
to be true  by unpublished work of 
Bartels and L\"uck~\cite{Bartels-Lueck(2008Borel)} and the
Bost Conjecture with $C^*$-coefficients has been proved by
Bartels-Echterhoff-L\"uck~\cite{Bartels-Echterhoff-Lueck(2007colim)}.

\end{enumerate}


\subsection{Quasiisometry}

\begin{enumerate}

\item Are there finitely presented groups that are quasiisometric
such that one is solvable but the other is not?
 See Remark~\ref{rem:virtually_solvable_groups}.

\item Is the property of being poly-cyclic invariant under quasiisometry?

\item Is the Mal'cev completion of a finitely generated torsionfree nilpotent
group an invariant under quasiisometry?
See Conjecture~\ref{con:folk_conjecture}.

\item Are the Novikov-Shubin invariants or the vanishing of the $L^2$-torsion
invariants under quasiisometry?
See Remark~\ref{rem:L2-invariants}.

\end{enumerate}


\typeout{-------------------- References
  -------------------------------}


\begin{thebibliography}{100}

\bibitem{Alibegovic+Bestvina(2006)}
E.~Alibegovi{\'c} and M.~Bestvina.
\newblock Limit groups are {$\rm CAT(0)$}.
\newblock {\em J. London Math. Soc. (2)}, 74(1):259--272, 2006.

\bibitem{Allcock+Gersten(1999)}
D.~J. Allcock and S.~M. Gersten.
\newblock A homological characterization of hyperbolic groups.
\newblock {\em Invent. Math.}, 135(3):723--742, 1999.

\bibitem{Alonso(1994)}
J.~M. Alonso.
\newblock Finiteness conditions on groups and quasi-isometries.
\newblock {\em J. Pure Appl. Algebra}, 95(2):121--129, 1994.

\bibitem{Bartels-Echterhoff-Lueck(2007colim)}
A.~Bartels, S.~Echterhoff, and W.~L\"uck.
\newblock Inheritance of isomorphism conjectures under colimits.
\newblock Preprintreihe SFB 478 --- Geometrische Strukturen in der Mathematik,
  Heft 452, M\"unster, arXiv:math.KT/0702460, to appear in the Proceedings of
  the conference on K-theory and non-commutative geometry in Valladolid,
  August/September 2006, 2007.

\bibitem{Bartels-Lueck(2008Borel)}
A.~Bartels and W.~L\"uck.
\newblock The {B}orel conjecture for hyperbolic and {C}{A}{T}(0)-groups.
\newblock in preparation, 2008.

\bibitem{Bartels-Lueck-Weinberger(2008)}
A.~Bartels, W.~L\"uck, and S.~Weinberger.
\newblock On hyperbolic groups with spheres as boundary.
\newblock in preparation, 2008.

\bibitem{Baum-Connes-Higson(1994)}
P.~Baum, A.~Connes, and N.~Higson.
\newblock Classifying space for proper actions and ${K}$-theory of group
  ${C}\sp \ast$-algebras.
\newblock In {\em $C\sp \ast$-algebras: 1943--1993 (San Antonio, TX, 1993)},
  pages 240--291. Amer. Math. Soc., Providence, RI, 1994.

\bibitem{Baumslag(1993)}
G.~Baumslag.
\newblock {\em Topics in combinatorial group theory}.
\newblock Birkh\"auser Verlag, Basel, 1993.

\bibitem{Behrstock-Kleiner-Minsky-Mosher(2008)}
J.~Behrstock, B.~Kleiner, Y.~Minsky, and L.~Mosher.
\newblock Geometry and rigidity of mapping class groups.
\newblock Preprint, arXiv:math.GT/0801.2006, 2008.

\bibitem{Bell(2005)}
R.~W. Bell.
\newblock Three-dimensional {FC} {A}rtin groups are {CAT}(0).
\newblock {\em Geom. Dedicata}, 113:21--53, 2005.

\bibitem{Benedetti-Petronio(1992)}
R.~Benedetti and C.~Petronio.
\newblock {\em Lectures on hyperbolic geometry}.
\newblock Springer-Verlag, Berlin, 1992.

\bibitem{Bestvina-Mess(1991)}
M.~Bestvina and G.~Mess.
\newblock The boundary of negatively curved groups.
\newblock {\em J. Amer. Math. Soc.}, 4(3):469--481, 1991.

\bibitem{Block+Weinberger(1997)}
J.~Block and S.~Weinberger.
\newblock Large scale homology theories and geometry.
\newblock In {\em Geometric topology (Athens, GA, 1993)}, volume~2 of {\em
  AMS/IP Stud. Adv. Math.}, pages 522--569. Amer. Math. Soc., Providence, RI,
  1997.

\bibitem{Bonk-Kleiner(2005)}
M.~Bonk and B.~Kleiner.
\newblock Conformal dimension and {G}romov hyperbolic groups with 2-sphere
  boundary.
\newblock {\em Geom. Topol.}, 9:219--246 (electronic), 2005.

\bibitem{Boudon(1995)}
M.~Bourdon.
\newblock Structure conforme au bord et flot g\'eod\'esique d'un {${\rm
  CAT}(-1)$}-espace.
\newblock {\em Enseign. Math. (2)}, 41(1-2):63--102, 1995.

\bibitem{Bowditch(1991)}
B.~H. Bowditch.
\newblock Notes on {G}romov's hyperbolicity criterion for path-metric spaces.
\newblock In {\em Group theory from a geometrical viewpoint (Trieste, 1990)},
  pages 64--167. World Sci. Publishing, River Edge, NJ, 1991.

\bibitem{Bowditch(1998)}
B.~H. Bowditch.
\newblock Cut points and canonical splittings of hyperbolic groups.
\newblock {\em Acta Math.}, 180(2):145--186, 1998.

\bibitem{Brick(1993)}
S.~G. Brick.
\newblock Quasi-isometries and ends of groups.
\newblock {\em J. Pure Appl. Algebra}, 86(1):23--33, 1993.

\bibitem{Bridson-Gersten(1996)}
M.~R. Bridson and S.~M. Gersten.
\newblock The optimal isoperimetric inequality for torus bundles over the
  circle.
\newblock {\em Quart. J. Math. Oxford Ser. (2)}, 47(185):1--23, 1996.

\bibitem{Bridson-Haefliger(1999)}
M.~R. Bridson and A.~Haefliger.
\newblock {\em Metric spaces of non-positive curvature}.
\newblock Springer-Verlag, Berlin, 1999.
\newblock Die Grundlehren der mathematischen Wissenschaften, Band 319.

\bibitem{Brown(1982)}
K.~S. Brown.
\newblock {\em Cohomology of groups}, volume~87 of {\em Graduate Texts in
  Mathematics}.
\newblock Springer-Verlag, New York, 1982.

\bibitem{Bumagina(2004)}
I.~Bumagina.
\newblock The {H}opf property for subgroups of hyperbolic groups.
\newblock {\em Geom. Dedicata}, 106:211--230, 2004.

\bibitem{Burger+Monod(2002)}
M.~Burger and N.~Monod.
\newblock Continuous bounded cohomology and applications to rigidity theory.
\newblock {\em Geom. Funct. Anal.}, 12(2):219--280, 2002.

\bibitem{Caprace-Monod(2008)}
P.-E. Caprace and N.~Monod.
\newblock Some properties of non-positively curved lattices.
\newblock Preprint, arXiv:math.GR/0806.0156v1, 2008.

\bibitem{Carlsson-Pedersen(1995a)}
G.~Carlsson and E.~K. Pedersen.
\newblock Controlled algebra and the {N}ovikov conjectures for ${K}$- and
  ${L}$-theory.
\newblock {\em Topology}, 34(3):731--758, 1995.

\bibitem{Casson-Jungreis(1994)}
A.~Casson and D.~Jungreis.
\newblock Convergence groups and {S}eifert fibered $3$-manifolds.
\newblock {\em Invent. Math.}, 118(3):441--456, 1994.

\bibitem{Charney-Davis(1995)}
R.~M. Charney and M.~W. Davis.
\newblock Strict hyperbolization.
\newblock {\em Topology}, 34(2):329--350, 1995.

\bibitem{Cohen(1989)}
D.~E. Cohen.
\newblock {\em Combinatorial group theory: a topological approach}.
\newblock Cambridge University Press, Cambridge, 1989.

\bibitem{Cowling+Haagerup(1989)}
M.~Cowling and U.~Haagerup.
\newblock Completely bounded multipliers of the {F}ourier algebra of a simple
  {L}ie group of real rank one.
\newblock {\em Invent. Math.}, 96(3):507--549, 1989.

\bibitem{Culler-Vogtmann(1986)}
M.~Culler and K.~Vogtmann.
\newblock Moduli of graphs and automorphisms of free groups.
\newblock {\em Invent. Math.}, 84(1):91--119, 1986.

\bibitem{Davis(2000Poin)}
M.~W. Davis.
\newblock Poincar\'e duality groups.
\newblock In {\em Surveys on surgery theory, Vol. 1}, volume 145 of {\em Ann.
  of Math. Stud.}, pages 167--193. Princeton Univ. Press, Princeton, NJ, 2000.

\bibitem{Davis(2008cox)}
M.~W. Davis.
\newblock {\em The geometry and topology of {C}oxeter groups}, volume~32 of
  {\em London Mathematical Society Monographs Series}.
\newblock Princeton University Press, Princeton, NJ, 2008.

\bibitem{Davis-Januszkiewicz(1991)}
M.~W. Davis and T.~Januszkiewicz.
\newblock Hyperbolization of polyhedra.
\newblock {\em J. Differential Geom.}, 34(2):347--388, 1991.

\bibitem{delaHarpe(2000)}
P.~de~la Harpe.
\newblock {\em Topics in geometric group theory}.
\newblock University of Chicago Press, Chicago, IL, 2000.

\bibitem{Dicks-Dunwoody(1989)}
W.~Dicks and M.~J. Dunwoody.
\newblock {\em Groups acting on graphs}.
\newblock Cambridge University Press, Cambridge, 1989.

\bibitem{Dioubina(1999)}
A.~Diubina.
\newblock Instability of the virtual solvability and the property of being
  virtually torsion-free for quasi-isometric groups.
\newblock Preprint, arXiv:math.GR/9911099, 1999.

\bibitem{Drutu(2002)}
C.~Dru{\c{t}}u.
\newblock Quasi-isometry invariants and asymptotic cones.
\newblock {\em Internat. J. Algebra Comput.}, 12(1-2):99--135, 2002.
\newblock International Conference on Geometric and Combinatorial Methods in
  Group Theory and Semigroup Theory (Lincoln, NE, 2000).

\bibitem{Dunwoody+Sageev(1999)}
M.~J. Dunwoody and M.~E. Sageev.
\newblock J{SJ}-splittings for finitely presented groups over slender groups.
\newblock {\em Invent. Math.}, 135(1):25--44, 1999.

\bibitem{Epstein(2008)}
I.~Epstein.
\newblock Orbit equivalent actions of non-amenable groups.
\newblock arXiv:math.GR/0707.4215, 2008.

\bibitem{Eskin(1998)}
A.~Eskin.
\newblock Quasi-isometric rigidity of nonuniform lattices in higher rank
  symmetric spaces.
\newblock {\em J. Amer. Math. Soc.}, 11(2):321--361, 1998.

\bibitem{Eskin+Farb(1998)}
A.~Eskin and B.~Farb.
\newblock Quasi-flats in {$H\sp 2\times H\sp 2$}.
\newblock In {\em Lie groups and ergodic theory (Mumbai, 1996)}, volume~14 of
  {\em Tata Inst. Fund. Res. Stud. Math.}, pages 75--103. Tata Inst. Fund.
  Res., Bombay, 1998.

\bibitem{Farb(1997)}
B.~Farb.
\newblock The quasi-isometry classification of lattices in semisimple {L}ie
  groups.
\newblock {\em Math. Res. Lett.}, 4(5):705--717, 1997.

\bibitem{Farb+Mosher(1999)}
B.~Farb and L.~Mosher.
\newblock Quasi-isometric rigidity for the solvable {B}aumslag-{S}olitar
  groups. {II}.
\newblock {\em Invent. Math.}, 137(3):613--649, 1999.

\bibitem{Farb+Mosher(2000)}
B.~Farb and L.~Mosher.
\newblock On the asymptotic geometry of abelian-by-cyclic groups.
\newblock {\em Acta Math.}, 184(2):145--202, 2000.

\bibitem{Farb+Schwartz(1996)}
B.~Farb and R.~Schwartz.
\newblock The large-scale geometry of {H}ilbert modular groups.
\newblock {\em J. Differential Geom.}, 44(3):435--478, 1996.

\bibitem{Farley(2003)}
D.~S. Farley.
\newblock Proper isometric actions of {T}hompson's groups on {H}ilbert space.
\newblock {\em Int. Math. Res. Not.}, 45:2409--2414, 2003.

\bibitem{Farrell-Jones(1993a)}
F.~T. Farrell and L.~E. Jones.
\newblock Isomorphism conjectures in algebraic ${K}$-theory.
\newblock {\em J. Amer. Math. Soc.}, 6(2):249--297, 1993.

\bibitem{Foelner(1955)}
E.~F{\o}lner.
\newblock On groups with full {B}anach mean value.
\newblock {\em Math. Scand.}, 3:243--254, 1955.

\bibitem{Freden(1995)}
E.~M. Freden.
\newblock Negatively curved groups have the convergence property. {I}.
\newblock {\em Ann. Acad. Sci. Fenn. Ser. A I Math.}, 20(2):333--348, 1995.

\bibitem{Furman(1999a)}
A.~Furman.
\newblock Gromov's measure equivalence and rigidity of higher rank lattices.
\newblock {\em Ann. of Math. (2)}, 150(3):1059--1081, 1999.

\bibitem{Furman(1999b)}
A.~Furman.
\newblock Orbit equivalence rigidity.
\newblock {\em Ann. of Math. (2)}, 150(3):1083--1108, 1999.

\bibitem{Gabai(1991)}
D.~Gabai.
\newblock Convergence groups are {F}uchsian groups.
\newblock {\em Bull. Amer. Math. Soc. (N.S.)}, 25(2):395--402, 1991.

\bibitem{Gaboriau(2002a)}
D.~Gaboriau.
\newblock Invariants {$l\sp 2$} de relations d'\'equivalence et de groupes.
\newblock {\em Publ. Math. Inst. Hautes \'Etudes Sci.}, (95):93--150, 2002.

\bibitem{Gaboriau(2002b)}
D.~Gaboriau.
\newblock On orbit equivalence of measure preserving actions.
\newblock In {\em Rigidity in dynamics and geometry (Cambridge, 2000)}, pages
  167--186. Springer, Berlin, 2002.

\bibitem{Gaboriau(2005)}
D.~Gaboriau.
\newblock Examples of groups that are measure equivalent to the free group.
\newblock {\em Ergodic Theory Dynam. Systems}, 25(6):1809--1827, 2005.

\bibitem{Gaboriau-Popa(2005)}
D.~Gaboriau and S.~Popa.
\newblock An uncountable family of nonorbit equivalent actions of {$\Bbb F\sb
  n$}.
\newblock {\em J. Amer. Math. Soc.}, 18(3):547--559 (electronic), 2005.

\bibitem{Gersten(1993)}
S.~M. Gersten.
\newblock Quasi-isometry invariance of cohomological dimension.
\newblock {\em C. R. Acad. Sci. Paris S\'er. I Math.}, 316(5):411--416, 1993.

\bibitem{Ghys-Harpe(1990)}
{\'E}.~Ghys and P.~de~la Harpe, editors.
\newblock {\em Sur les groupes hyperboliques d'apr\`es {M}ikhael {G}romov}.
\newblock Birkh\"auser Boston Inc., Boston, MA, 1990.
\newblock Papers from the Swiss Seminar on Hyperbolic Groups held in Bern,
  1988.

\bibitem{Ghys+delaHarpe(1991)}
{\'E}.~Ghys and P.~de~la Harpe.
\newblock Infinite groups as geometric objects (after {G}romov).
\newblock In {\em Ergodic theory, symbolic dynamics, and hyperbolic spaces
  (Trieste, 1989)}, Oxford Sci. Publ., pages 299--314. Oxford Univ. Press, New
  York, 1991.

\bibitem{Gromov(1981c)}
M.~Gromov.
\newblock Groups of polynomial growth and expanding maps.
\newblock {\em Inst. Hautes \'Etudes Sci. Publ. Math.}, 53:53--73, 1981.

\bibitem{Gromov(1981)}
M.~Gromov.
\newblock {\em Structures m\'etriques pour les vari\'et\'es riemanniennes}.
\newblock CEDIC, Paris, 1981.
\newblock Edited by J. Lafontaine and P. Pansu.

\bibitem{Gromov(1987)}
M.~Gromov.
\newblock Hyperbolic groups.
\newblock In {\em Essays in group theory}, pages 75--263. Springer-Verlag, New
  York, 1987.

\bibitem{Gromov(1993)}
M.~Gromov.
\newblock Asymptotic invariants of infinite groups.
\newblock In {\em Geometric group theory, Vol.\ 2 (Sussex, 1991)}, pages
  1--295. Cambridge Univ. Press, Cambridge, 1993.

\bibitem{Gromov(2003)}
M.~Gromov.
\newblock Random walk in random groups.
\newblock {\em Geom. Funct. Anal.}, 13(1):73--146, 2003.

\bibitem{Hamenstaedt(2007)}
U.~Hamenst\"adt.
\newblock Geometry of the mapping class groups iii: Quasi-isometry rigidity.
\newblock Preprint, arXiv:math.GT/0512429v2, 2007.

\bibitem{Higson-Kasparov(2001)}
N.~Higson and G.~Kasparov.
\newblock ${E}$-theory and ${K}{K}$-theory for groups which act properly and
  isometrically on {H}ilbert space.
\newblock {\em Invent. Math.}, 144(1):23--74, 2001.

\bibitem{Higson-Lafforgue-Skandalis(2002)}
N.~Higson, V.~Lafforgue, and G.~Skandalis.
\newblock Counterexamples to the {B}aum-{C}onnes conjecture.
\newblock {\em Geom. Funct. Anal.}, 12(2):330--354, 2002.

\bibitem{Hjorth(2005)}
G.~Hjorth.
\newblock A converse to {D}ye's theorem.
\newblock {\em Trans. Amer. Math. Soc.}, 357(8):3083--3103 (electronic), 2005.

\bibitem{Kapovich+Benakli(2002)}
I.~Kapovich and N.~Benakli.
\newblock Boundaries of hyperbolic groups.
\newblock In {\em Combinatorial and geometric group theory (New York,
  2000/Hoboken, NJ, 2001)}, volume 296 of {\em Contemp. Math.}, pages 39--93.
  Amer. Math. Soc., Providence, RI, 2002.

\bibitem{Kapovich+Kleiner(2007)}
M.~Kapovich and B.~Kleiner.
\newblock The weak hyperbolization conjecture for 3-dimensional
  {CAT}(0)-groups.
\newblock {\em Groups Geom. Dyn.}, 1(1):61--79, 2007.

\bibitem{Kapovich-Leeb(1995)}
M.~Kapovich and B.~Leeb.
\newblock On asymptotic cones and quasi-isometry classes of fundamental groups
  of $3$-manifolds.
\newblock {\em Geom. Funct. Anal.}, 5(3):582--603, 1995.

\bibitem{Kleiner(2007)}
B.~Kleiner.
\newblock A new proof of {G}romov's theorem on groups of polynomial growth.
\newblock Preprint, arXiv:math.GR/ 0710.4593v4, 2007.

\bibitem{Kleiner-Leeb(1997)}
B.~Kleiner and B.~Leeb.
\newblock Rigidity of quasi-isometries for symmetric spaces and {E}uclidean
  buildings.
\newblock {\em Inst. Hautes \'Etudes Sci. Publ. Math.}, 86:115--197 (1998),
  1997.

\bibitem{Kropholler-Mislin(1998)}
P.~H. Kropholler and G.~Mislin.
\newblock Groups acting on finite-dimensional spaces with finite stabilizers.
\newblock {\em Comment. Math. Helv.}, 73(1):122--136, 1998.

\bibitem{Lafforgue(1999)}
V.~Lafforgue.
\newblock Compl\'ements \`a la d\'emonstration de la conjecture de
  {B}aum-{C}onnes pour certains groupes poss\'edant la propri\'et\'e ({T}).
\newblock {\em C. R. Acad. Sci. Paris S\'er. I Math.}, 328(3):203--208, 1999.

\bibitem{Lafforgue(2001)}
V.~Lafforgue.
\newblock Banach {$KK$}-theory and the {B}aum-{C}onnes conjecture.
\newblock In {\em European Congress of Mathematics, Vol. II (Barcelona, 2000)},
  volume 202 of {\em Progr. Math.}, pages 31--46. Birkh\"auser, Basel, 2001.

\bibitem{Leary-Nucinkis(2001a)}
I.~J. Leary and B.~E.~A. Nucinkis.
\newblock Every {CW}-complex is a classifying space for proper bundles.
\newblock {\em Topology}, 40(3):539--550, 2001.

\bibitem{Leary-Nucinkis(2003)}
I.~J. Leary and B.~E.~A. Nucinkis.
\newblock Some groups of type {$VF$}.
\newblock {\em Invent. Math.}, 151(1):135--165, 2003.

\bibitem{Lueck(1989)}
W.~L{\"u}ck.
\newblock {\em Transformation groups and algebraic ${K}$-theory}.
\newblock Springer-Verlag, Berlin, 1989.
\newblock Mathematica Gottingensis.

\bibitem{Lueck(2000a)}
W.~L{\"u}ck.
\newblock The type of the classifying space for a family of subgroups.
\newblock {\em J. Pure Appl. Algebra}, 149(2):177--203, 2000.

\bibitem{Lueck(2002b)}
W.~L{\"u}ck.
\newblock Chern characters for proper equivariant homology theories and
  applications to ${K}$- and ${L}$-theory.
\newblock {\em J. Reine Angew. Math.}, 543:193--234, 2002.

\bibitem{Lueck(2002)}
W.~L{\"u}ck.
\newblock {\em {$L\sp 2$}-invariants: theory and applications to geometry and
  {$K$}-theory}, volume~44 of {\em Ergebnisse der Mathematik und ihrer
  Grenzgebiete. 3. Folge. A Series of Modern Surveys in Mathematics [Results in
  Mathematics and Related Areas. 3rd Series. A Series of Modern Surveys in
  Mathematics]}.
\newblock Springer-Verlag, Berlin, 2002.

\bibitem{Lueck(2002d)}
W.~L{\"u}ck.
\newblock The relation between the {B}aum-{C}onnes conjecture and the trace
  conjecture.
\newblock {\em Invent. Math.}, 149(1):123--152, 2002.

\bibitem{Lueck(2005heis)}
W.~L{\"u}ck.
\newblock {$K$}- and {$L$}-theory of the semi-direct product of the discrete
  3-dimensional {H}eisenberg group by {${\Bbb Z}/4$}.
\newblock {\em Geom. Topol.}, 9:1639--1676 (electronic), 2005.

\bibitem{Lueck(2005s)}
W.~L{\"u}ck.
\newblock Survey on classifying spaces for families of subgroups.
\newblock In {\em Infinite groups: geometric, combinatorial and dynamical
  aspects}, volume 248 of {\em Progr. Math.}, pages 269--322. Birkh\"auser,
  Basel, 2005.

\bibitem{Lueck(2007)}
W.~L{\"u}ck.
\newblock Rational computations of the topological {$K$}-theory of classifying
  spaces of discrete groups.
\newblock {\em J. Reine Angew. Math.}, 611:163--187, 2007.

\bibitem{Lueck(2009catevcyc)}
W.~L{\"u}ck.
\newblock On the classifying space of the family of virtually cyclic subgroups
  for {C}{A}{T}(0)-groups.
\newblock in preparation, 2009.

\bibitem{Lueck-Meintrup(2000)}
W.~L{\"u}ck and D.~Meintrup.
\newblock On the universal space for group actions with compact isotropy.
\newblock In {\em Geometry and topology: Aarhus (1998)}, pages 293--305. Amer.
  Math. Soc., Providence, RI, 2000.

\bibitem{Lueck-Reich(2005)}
W.~L{\"u}ck and H.~Reich.
\newblock The {B}aum-{C}onnes and the {F}arrell-{J}ones conjectures in {$K$}-
  and {$L$}-theory.
\newblock In {\em Handbook of $K$-theory. Vol. 1, 2}, pages 703--842. Springer,
  Berlin, 2005.

\bibitem{Lueck-Stamm(2000)}
W.~L{\"u}ck and R.~Stamm.
\newblock Computations of ${K}$- and ${L}$-theory of cocompact planar groups.
\newblock {\em $K$-Theory}, 21(3):249--292, 2000.

\bibitem{Lyndon-Schupp(1977)}
R.~C. Lyndon and P.~E. Schupp.
\newblock {\em Combinatorial group theory}.
\newblock Springer-Verlag, Berlin, 1977.
\newblock Ergebnisse der Mathematik und ihrer Grenzgebiete, Band 89.

\bibitem{Meintrup-Schick(2002)}
D.~Meintrup and T.~Schick.
\newblock A model for the universal space for proper actions of a hyperbolic
  group.
\newblock {\em New York J. Math.}, 8:1--7 (electronic), 2002.

\bibitem{Mineyev(2002boundcoh)}
I.~Mineyev.
\newblock Bounded cohomology characterizes hyperbolic groups.
\newblock {\em Q. J. Math.}, 53(1):59--73, 2002.

\bibitem{Monod+Shalom(2006)}
N.~Monod and Y.~Shalom.
\newblock Orbit equivalence rigidity and bounded cohomology.
\newblock {\em Ann. of Math. (2)}, 164(3):825--878, 2006.

\bibitem{Moussong(1987)}
G.~Moussong.
\newblock {\em Hyperbolic {C}oxeter groups}.
\newblock Ph. d. thesis, The Ohio State University, 1987.

\bibitem{Nomizu(1954)}
K.~Nomizu.
\newblock On the cohomology of compact homogeneous spaces of nilpotent {L}ie
  groups.
\newblock {\em Ann. of Math. (2)}, 59:531--538, 1954.

\bibitem{Olshanskii-Osin-Sapir(2007)}
A.~Y. Ol'shankskii, D.~Osin, and M.~Sapir.
\newblock Lacunary hyperbolic groups.
\newblock arXiv:math.GR/0701365v1, 2007.

\bibitem{Olshanskii(1992)}
A.~Y. Ol'shanskii.
\newblock Almost every group is hyperbolic.
\newblock {\em Internat. J. Algebra Comput.}, 2(1):1--17, 1992.

\bibitem{Olshanskii(1993)}
A.~Y. Ol'shanski{\u\i}.
\newblock On residualing homomorphisms and {$G$}-subgroups of hyperbolic
  groups.
\newblock {\em Internat. J. Algebra Comput.}, 3(4):365--409, 1993.

\bibitem{Olshanskii+Sapir(2006)}
A.~Y. Olshanskii and M.~V. Sapir.
\newblock Groups with non-simply connected asymptotic cones.
\newblock In {\em Topological and asymptotic aspects of group theory}, volume
  394 of {\em Contemp. Math.}, pages 203--208. Amer. Math. Soc., Providence,
  RI, 2006.

\bibitem{Ornstein-Weiss(1980)}
D.~S. Ornstein and B.~Weiss.
\newblock Ergodic theory of amenable group actions. {I}. {T}he {R}ohlin lemma.
\newblock {\em Bull. Amer. Math. Soc. (N.S.)}, 2(1):161--164, 1980.

\bibitem{Ozawa(2007)}
N.~Ozawa.
\newblock Weak amenability of hyperbolic groups.
\newblock Preprint, arXiv:math.FA/0704.1635, 2007.

\bibitem{Ozawa+Popa(2004)}
N.~Ozawa and S.~Popa.
\newblock Some prime factorization results for type {${\rm II}\sb 1$} factors.
\newblock {\em Invent. Math.}, 156(2):223--234, 2004.

\bibitem{Pansu(1989c)}
P.~Pansu.
\newblock M\'etriques de {C}arnot-{C}arath\'eodory et quasiisom\'etries des
  espaces sym\'etriques de rang un.
\newblock {\em Ann. of Math. (2)}, 129(1):1--60, 1989.

\bibitem{Pansu(1995)}
P.~Pansu.
\newblock Cohomologie ${L^p}$: Invariance sous quasiisometries.
\newblock Preprint, Orsay, 1995.

\bibitem{Papasoglu(1995trees)}
P.~Papasoglu.
\newblock Homogeneous trees are bi-{L}ipschitz equivalent.
\newblock {\em Geom. Dedicata}, 54(3):301--306, 1995.

\bibitem{Papasoglu(2005)}
P.~Papasoglu.
\newblock Quasi-isometry invariance of group splittings.
\newblock {\em Ann. of Math. (2)}, 161(2):759--830, 2005.

\bibitem{Papasoglu+Whyte(2002)}
P.~Papasoglu and K.~Whyte.
\newblock Quasi-isometries between groups with infinitely many ends.
\newblock {\em Comment. Math. Helv.}, 77(1):133--144, 2002.

\bibitem{Paterson(1988)}
A.~L.~T. Paterson.
\newblock {\em Amenability}.
\newblock American Mathematical Society, Providence, RI, 1988.

\bibitem{Paulin(1996)}
F.~Paulin.
\newblock Un groupe hyperbolique est d\'etermin\'e par son bord.
\newblock {\em J. London Math. Soc. (2)}, 54(1):50--74, 1996.

\bibitem{Popa(2004)}
S.~Popa.
\newblock On the fundamental group of type {$\rm II\sb 1$} factors.
\newblock {\em Proc. Natl. Acad. Sci. USA}, 101(3):723--726 (electronic), 2004.

\bibitem{Popa(2006_betti)}
S.~Popa.
\newblock On a class of type {${\rm II}\sb 1$} factors with {B}etti numbers
  invariants.
\newblock {\em Ann. of Math. (2)}, 163(3):809--899, 2006.

\bibitem{Popa(2006_rigidity_bernoulli)}
S.~Popa.
\newblock Some rigidity results for non-commutative {B}ernoulli shifts.
\newblock {\em J. Funct. Anal.}, 230(2):273--328, 2006.

\bibitem{Popa(2006_strong_rigidity_actions_I)}
S.~Popa.
\newblock Strong rigidity of {$\rm II\sb 1$} factors arising from malleable
  actions of {$w$}-rigid groups. {I}.
\newblock {\em Invent. Math.}, 165(2):369--408, 2006.

\bibitem{Popa(2006_strong_rigidity_actions_II}
S.~Popa.
\newblock Strong rigidity of {$\rm II\sb 1$} factors arising from malleable
  actions of {$w$}-rigid groups. {II}.
\newblock {\em Invent. Math.}, 165(2):409--451, 2006.

\bibitem{Popa(2007_rigidity)}
S.~Popa.
\newblock Cocycle and orbit equivalence superrigidity for malleable actions of
  {$w$}-rigid groups.
\newblock {\em Invent. Math.}, 170(2):243--295, 2007.

\bibitem{Popa(2007)}
S.~Popa.
\newblock Deformation and rigidity for group actions and von {N}eumann
  algebras.
\newblock In {\em International Congress of Mathematicians. Vol. I}, pages
  445--477. Eur. Math. Soc., Z\"urich, 2007.

\bibitem{Popa+Vaes(2008)}
S.~Popa and S.~Vaes.
\newblock Strong rigidity of generalized {B}ernoulli actions and computations
  of their symmetry groups.
\newblock {\em Adv. Math.}, 217(2):833--872, 2008.

\bibitem{Rips-Sela(1997)}
E.~Rips and Z.~Sela.
\newblock Cyclic splittings of finitely presented groups and the canonical
  {J}{S}{J} decomposition.
\newblock {\em Ann. of Math. (2)}, 146(1):53--109, 1997.

\bibitem{Roe(2005asym)}
J.~Roe.
\newblock Hyperbolic groups have finite asymptotic dimension.
\newblock {\em Proc. Amer. Math. Soc.}, 133(9):2489--2490 (electronic), 2005.

\bibitem{Rosenthal(2004)}
D.~Rosenthal.
\newblock Splitting with continuous control in algebraic {$K$}-theory.
\newblock {\em $K$-Theory}, 32(2):139--166, 2004.

\bibitem{Rosenthal-Schuetz(2005)}
D.~Rosenthal and D.~Sch{\"u}tz.
\newblock On the algebraic {$K$}- and {$L$}-theory of word hyperbolic groups.
\newblock {\em Math. Ann.}, 332(3):523--532, 2005.

\bibitem{Sauer(2006)}
R.~Sauer.
\newblock Homological invariants and quasi-isometry.
\newblock {\em Geom. Funct. Anal.}, 16(2):476--515, 2006.

\bibitem{Schwartz(1995)}
R.~E. Schwartz.
\newblock The quasi-isometry classification of rank one lattices.
\newblock {\em Inst. Hautes \'Etudes Sci. Publ. Math.}, 82:133--168 (1996),
  1995.

\bibitem{Schwartz(1996)}
R.~E. Schwartz.
\newblock Quasi-isometric rigidity and {D}iophantine approximation.
\newblock {\em Acta Math.}, 177(1):75--112, 1996.

\bibitem{Scott+Swarup(2003)}
P.~Scott and G.~A. Swarup.
\newblock Regular neighbourhoods and canonical decompositions for groups.
\newblock {\em Ast\'erisque}, 289:vi+233, 2003.

\bibitem{Sela(1995)}
Z.~Sela.
\newblock The isomorphism problem for hyperbolic groups. {I}.
\newblock {\em Ann. of Math. (2)}, 141(2):217--283, 1995.

\bibitem{Sela(1997)}
Z.~Sela.
\newblock Structure and rigidity in ({G}romov) hyperbolic groups and discrete
  groups in rank {$1$} {L}ie groups. {II}.
\newblock {\em Geom. Funct. Anal.}, 7(3):561--593, 1997.

\bibitem{Serre(1971)}
J.-P. Serre.
\newblock Cohomologie des groupes discrets.
\newblock In {\em Prospects in mathematics (Proc. Sympos., Princeton Univ.,
  Princeton, N.J., 1970)}, pages 77--169. Ann. of Math. Studies, No. 70.
  Princeton Univ. Press, Princeton, N.J., 1971.

\bibitem{Serre(1980)}
J.-P. Serre.
\newblock {\em Trees}.
\newblock Springer-Verlag, Berlin, 1980.
\newblock Translated from the French by J.~Stillwell.

\bibitem{Shalom(2004)}
Y.~Shalom.
\newblock Harmonic analysis, cohomology, and the large-scale geometry of
  amenable groups.
\newblock {\em Acta Math.}, 192(2):119--185, 2004.

\bibitem{Shalom(2005)}
Y.~Shalom.
\newblock Measurable group theory.
\newblock In {\em European Congress of Mathematics}, pages 391--423. Eur. Math.
  Soc., Z\"urich, 2005.

\bibitem{Thomasson+Woess(1993)}
C.~Thomassen and W.~Woess.
\newblock Vertex-transitive graphs and accessibility.
\newblock {\em J. Combin. Theory Ser. B}, 58(2):248--268, 1993.

\bibitem{Dieck(1972)}
T.~tom Dieck.
\newblock Orbittypen und \"aquivariante {H}omologie. {I}.
\newblock {\em Arch. Math. (Basel)}, 23:307--317, 1972.

\bibitem{Dieck(1987)}
T.~tom Dieck.
\newblock {\em Transformation groups}.
\newblock Walter de Gruyter \& Co., Berlin, 1987.

\bibitem{Dries+Wilkie(1984)}
L.~van~den Dries and A.~J. Wilkie.
\newblock Gromov's theorem on groups of polynomial growth and elementary logic.
\newblock {\em J. Algebra}, 89(2):349--374, 1984.

\bibitem{Vogtmann(2002)}
K.~Vogtmann.
\newblock Automorphisms of free groups and outer space.
\newblock In {\em Proceedings of the Conference on Geometric and Combinatorial
  Group Theory, Part I (Haifa, 2000)}, volume~94, pages 1--31, 2002.

\bibitem{Wegner(2001)}
C.~Wegner.
\newblock ${L^2}$-invariants of finite aspherical ${CW}$-complexes.
\newblock Preprintreihe SFB 478 --- Geometrische Strukturen in der Mathematik,
  Heft 152, M\"unster, 2001.

\bibitem{Whyte(1999)}
K.~Whyte.
\newblock Amenability, bi-{L}ipschitz equivalence, and the von {N}eumann
  conjecture.
\newblock {\em Duke Math. J.}, 99(1):93--112, 1999.

\bibitem{Zuk(2003)}
A.~{\.Z}uk.
\newblock Property ({T}) and {K}azhdan constants for discrete groups.
\newblock {\em Geom. Funct. Anal.}, 13(3):643--670, 2003.

\end{thebibliography}

\end{document}